\documentclass[reqno,11pt]{amsart}
\usepackage{latexsym}
\usepackage{epsfig}
\calclayout
\def\RE{{\mathbb R}}
\def\Om{{\Omega}}
\def\spb{\mbox{spec}_b}
\def\ci{C^{\infty}}
\def\be{\overline{e}}
\def\bI{\overline{I}}
\def\CX{{\mathbb C}}
\def\IZ{{\mathbb Z}}
\def\db{\overline{\partial}}
\def\p{\partial}
\def\dir{\partial\!\!\!/}
\renewcommand{\Im}{\mbox{Im}\,}
\renewcommand{\Re}{\mbox{Re}\,}
\newcommand{\Ga}{\Gamma}
\newcommand{\bX}{\overline{X}}
\newcommand{\bu}{\overline{u}}
\newcommand{\bS}{\overline{S}}
\newcommand{\hS}{\widehat{S}}
\newcommand{\hX}{\widehat{X}}
\newcommand{\bg}{\overline{g}}
\newcommand{\wg}{\widetilde{g}}
\newcommand{\hg}{\widehat{g}}

\newcommand{\CI}{\mathcal{I}}

\newcommand{\CO}{\mathcal{O}}
\newcommand{\bc}{\overline{c}}
\newcommand{\hc}{\widehat{c}}
\newcommand{\bxi}{\overline{\xi}}
\newcommand{\hxi}{\widehat{\xi}}

\newcommand{\bpsi}{\overline{\psi}}
\newcommand{\hpsi}{\widehat{\psi}}
\newcommand{\bh}{\overline{h}}

\newcommand{\bdel}{\overline{\nabla}}
\newcommand{\bL}{\overline{L}}
\newcommand{\hL}{\widehat{L}}

\newcommand{\bD}{\overline{D}}
\newcommand{\hD}{\widehat{D}}
\newcommand{\DD}{{\mathcal{D}}}

\newcommand{\ve}{\varepsilon}

\newcommand{\OCP}{\overline{\CX P^2}}
\newcommand{\tr}{\mbox{tr }}

\newtheorem{result}{Theorem}

\newtheorem{lem}[equation]{Lemma}
\newtheorem{lemma}[equation]{Lemma}
\newtheorem{propn}[equation]{Proposition}

\newtheorem{thm}[equation]{Theorem}

\newtheorem{conjecture}[equation]{Conjecture}

\theoremstyle{remark}
\newtheorem*{remark}{Remark}

\let\epsilon=\varepsilon

\DeclareMathOperator{\coker}{coker}

\begin{document}

\title[Complete anti-self-dual spaces]{Gluing theorems for complete
anti-self-dual spaces} 
\author{Alexei Kovalev and Michael Singer}
\address{Department of Mathematics and Statistics,
  University of Edinburgh, King's Buildings, Edinburgh EH9 3JZ,
  Great Britain}
\email{agk@maths.ed.ac.uk}
\email{michael@maths.ed.ac.uk}
\date{September 13, 2000}
\maketitle

\section{Introduction}

\subsection{Summary}
One of the special features of $4$-dimensional differential geometry is 
the existence of objects with self-dual (SD) or anti-self-dual (ASD) 
curvature. The objects in question can be connections in an auxiliary 
bundle over a $4$-manifold, leading to the study of instantons in 
Yang--Mills theory \cite{DK}, or as in this paper,
Riemannian metrics or conformal structures. Although such {\em ASD conformal 
structures} give absolute minima of the functional $c \mapsto 
\|W(c)\|^2_2$, where $W(c)$ denotes the Weyl tensor of the conformal 
structure $c$, variational methods are not well suited to the study of this 
problem, essentially because of its conformal invariance.  For this 
reason, {\em gluing theorems} provide a very important source of 
information about ASD conformal structures. Our purpose in this paper 
is to give some new and rather general gluing theorems for ASD and 
Hermitian--ASD conformal structures, following the method suggested by 
Floer in \cite{floer-conformal}.  The prototypical gluing theorem 
takes a pair $(X_j,c_j)$ ($j=1,2$) of compact conformally ASD 
$4$-manifolds and analyzes the problem of finding an ASD conformal 
structure $c$ on $X = X_1\sharp X_2$ that is `close to' $c_j$ in
suitable subsets $X_j\backslash B_j\subset X_1\sharp X_2$. In this 
situation there exist finite-dimensional vector spaces (the {\em 
obstruction spaces})  $H^2_{c_j}(X_j)$ whose vanishing is sufficient 
to guarantee the existence of $c$ with the desired properties. 
(If $H^2_{c_j}(X_j)\not=0$, then the gluing theorem yields a map from 
another finite-dimensional vector space into 
$H^2_{c_1}(X_1)\oplus H^2_{c_2}(X_2)$, the zeroes of which yield ASD 
conformal structures on $X_1\sharp X_2$.) 

The result just stated (gluing for compact conformally ASD spaces) 
was proved by Donaldson and Friedman \cite{DF} and in a very special 
case by Floer \cite{floer-conformal}.  The approach of \cite{DF} was 
to exploit the twistor description \cite{P,ahs} of conformally ASD 
spaces, to translate the gluing problem into one of deformation 
theory of complex singular spaces. Floer, on the other hand, worked 
directly with the $4$-manifolds and used some tools from the theory 
of elliptic operators on non-compact manifolds with cylindrical ends.

One of the motivations for the present work was the desire to extend 
the basic gluing theorem to handle the case where the $(X_j,c_j)$ are 
conformally ASD orbifolds with isolated singular points. Such an 
extension opens up the possibility of obtaining new examples of 
conformally ASD spaces by desingularizing conformally ASD orbifolds 
within the ASD category. Indeed,
the process of resolution of orbifold singularities amounts to taking 
a (generalized) connected sum with a suitable standard orbifold with 
precisely one singular point, and examples of such standard orbifolds 
are now known for many of the finite subgroups of $SO_4$ 
\cite{Kro,lebrun,GL}.  The simplest 
possible case of $\IZ_2$-singularities was studied in
\cite{lebrun-singer}, where the methods of \cite{DF} were extended to
give gluing theorems for these orbifolds. The method was further
extended to cover cyclic singularities, by Jian Zhou \cite{jz}, but
becomes increasingly complicated owing to the singularities developed
by the corresponding twistor spaces.

On general grounds, however, one should expect that Floer's 
analytical approach might provide a simpler framework for such 
generalized gluing theorems. Indeed, in that approach the first step 
is to blow up $X_j$ at the marked point $0_j$ at which the gluing takes 
place.  This blow-up results in a manifold with an {\em infinite
cylindrical end}, or more-or-less equivalently, a compact manifold with
boundary, equipped with a $b$-metric \cite{tapsit}. The cross-section
of the cylinder is diffeomorphic to the added boundary component, both
being the link in  $X_j$ of $0_j$, and  the singularity has disappeared
completely. Having reached this point, it is quite 
reasonable to take $4$-manifolds with {\em boundary}, equipped 
with conformally ASD $b$-metrics as the basic entities to glue, 
regarding compact manifolds and orbifolds as special cases.  Our main 
result here is indeed a gluing theorem for pairs of conformally ASD 
$b$-manifolds $(X_j,c_j)$, the connected sum being replaced by the 
`join' $X=X_1\cup_Y X_2$ of $X_1$ and $X_2$ across $Y \subset \p X_j$, 
where $Y$ is a common piece (union of connected components) of the 
boundaries of the $X_j$.  The precise statement involves a 
considerable notational overhead in the definition of the obstruction 
spaces and is deferred to \S\ref{mainthm}; suffice it to say that once the 
obstruction spaces have been correctly defined, the 
result is precisely analogous to the prototype mentioned before.  It 
should perhaps be emphasised that this theorem genuinely operates in 
the $b$-category, in that $\p X$ can be non-empty, in 
which case the gluing theorem produces a conformally ASD $b$-metric 
on $X$ or, in more traditional language, a complete conformally ASD 
metric on $X\backslash \p X$ with cylindrical asymptotics.

Such a gluing theorem immediately poses questions about the existence 
of conformally ASD $b$-metrics.  The first observation is that near 
each component $Y$  of the boundary, the metric must be asymptotic to
a conformally ASD product metric on $Y\times \RE$. Thus \cite{besse} 
$Y$ must have constant sectional curvature, so the new possibilities 
(not arising from orbifolds) are $Y =$  the 3-torus $T^3$ 
or $Y =$ a hyperbolic 
$3$-manifold. It seems that not much is known about the existence of 
conformally ASD $b$-metrics on manifolds with such boundary components,
so we argue in \S\ref{nct} that Taubes's method \cite{taubes-conformal,
taubes-glueing} can be adapted to 
yield conformally ASD $b$-metrics on $X\sharp N\OCP$ for large 
enough $N$, if $(X,g)$ is any Riemannian $b$-manifold such that 
$g$ is conformally flat near $\p X$.  In a companion paper 
\cite{KS}, we also give a simple example of a conformally ASD (in 
fact hyperK\"ahler) $b$-metric on a manifold with boundary equal to 
$T^3$.

In this paper we also study hermitian-ASD conformal structures on 
complex surfaces. These are particularly interesting because of their 
relation to scalar-flat K\"ahler geometry \cite{exsdm,LebS,KLP,asdk}. Such
K\"ahler metrics with zero scalar curvature are of 
interest from the point of view of Calabi's extremal metric 
programme \cite{besse}, and even for 
complex surfaces there is no systematic 
existence theory.  It is fortunate, then, that for surfaces,
scalar-flat K\"ahler metrics can be approached through hermitian-ASD 
conformal structures, and hence through gluing theorems.  As for the 
full ASD equations, we consider the general gluing problem for 
conformally ASD hermitian $b$-metrics on compact complex surfaces with 
boundary, and obtain similar results.  We also illustrate our general results with
a simple 
application, showing that the blow-up of $\CX^2$ at an arbitrary set 
of points $p_j$ admits scalar-flat K\"ahler metrics that are asymptotic to 
the Euclidean metric at $\infty$.  This generalizes LeBrun's explicit 
construction \cite[Theorem 1]{exsdm} of $S^1$-invariant
scalar-flat K\"ahler metrics on this blow-up when the $p_j$ lie on
a complex line in $\CX^2$.

It should perhaps be remarked that the ASD condition requires that an 
orientation be chosen on the underlying $4$-manifold, and accordingly 
some care has to be taken with the construction of $X_1\cup_Y X_2$ 
to ensure that this has an orientation that is compatible with the 
given orientations on the $X_j$. This is particularly true for gluing 
complex surfaces: a moment's reflection will convince 
the reader that the connected sum of 2 complex surfaces {\em  never} has
a complex structure compatible with the given complex structures on 
the summands. One will, however, be able to glue 
the asymptotic region of a suitable non-compact surface to the 
complement of a neighbourhood of a point in a compact surface. In 
this way one gets gluing theorems that give sufficient conditions for 
the blow-up of a compact scalar-flat K\"ahler metric to admit a 
scalar-flat K\"ahler metric or more generally results 
about resolution of 
orbifold singularities of scalar-flat K\"ahler metrics within the 
scalar-flat K\"ahler category.

\subsection{Strategy} Gluing theorems have been pursued vigorously in
many different 
contexts over the last 10 or so years, and there are many approaches
to the problem. For nonlinear problems like the ones studied in this
paper, involving the construction of connections or metrics
with prescribed
curvature, the strategy is always the same: construct a family of `approximate
solutions' on the join of the two spaces, and then use some variant of
the implicit function theorem to obtain a nearby genuine solution.  In
the notation of the rest of this paper, this family of solutions will
depend on a large parameter $\rho$, essentially the length of the neck
joining the two spaces. As $\rho \to \infty$, the approximate 
solution gets better and better.  In particular, the implicit function
theorem needs to be applied when $\rho$ is very large, and for this
some good control of the linearization of the problem is needed in
this limit.  For this reason a good deal of the work concerns the
behaviour of linear operators on manifolds with long necks. These
linear problems are of interest in their own right in the context of 
gluing formulae in index theory and for $\eta$-invariants: recent work
in this direction, from a point of view that is close to that of this
paper, can be found for example in \cite{mm,hmm}
For a more leisurely description of gluing theorems of
this type, the reader is referred to \cite[Chapter 7]{DK},
\cite{floer-conformal} or to \cite{taubes-glueing} for a recent survey
of gluing problems for instantons and ASD conformal structures. Some other
geometric (nonlinear)  gluing problems are surveyed in \cite{gglue}.

Having given these pointers to the literature, we shall
concentrate in the rest of this paper  on the technical details of
our particular problems, without too much further motivation.

\subsection{Contents}

In more detail, the remaining sections of the paper are as follows:
\begin{quote} 

\S2: The $b$-category is introduced, Fredholm properties of
$b$-differential operators are described, and gluing of $b$-metrics is
explained in detail.

\S3: Conformal geometry is recalled, with special reference to the ASD
equations in $4$ dimensions. The relevant PDE aspects of these
equations are described.

\S4: Linear aspects of our problem are discussed here, including a
rather general account of the behaviour of the kernel and cokernel of
elliptic operators on manifolds with long necks. The application of
these general results to the ASD problem is also treated along with
`comparison theorems' which allow us to compare the linearization of
the ASD equations on a compact orbifold with the linearization on the
corresponding `blown-up' $b$-manifold. 

\S5: Nonlinear aspects of the analysis appear here, centring around
the application of the implicit function theorem to obtain a weak
solution of the ASD equations, and elliptic regularity arguments to
show that this solution is $\ci$ (and has optimal behaviour at the
boundary, if there is one). 

\S6:  The main gluing theorems are summarized here, both for 
conformally ASD and conformally hermitian-ASD metrics. The 
construction of scalar-flat K\"ahler metrics on an arbitrary blow-up 
of $\CX^2$ also appears here.

\S7: The $b$-version of Taubes's existence theorem is given here.

\S8: A number of vanishing theorems for the
obstruction space $H^2_{c}(X)$ are collected here. 
\end{quote}

\subsection{Acknowledgements} 

We acknowledge useful conversations and encouragement from David
Calderbank, Dominic Joyce, Claude LeBrun, Rafe Mazzeo, 
Richard Melrose and Mario Micallef. Both authors were supported by
EPSRC while most of this work was carried out. 

Finally our approach owes much to the work of Andreas Floer
\cite{floer-conformal}.  We hope to honour this exceptional 
mathematician with the present work, by bringing his 
insights about gluing ASD conformal
structures to a wider audience.

\setcounter{equation}{0}
\section{Gluing $b$-manifolds}

In this section we recall Richard Melrose's approach to spaces with
cylindrical ends, and explain the relevance of these ideas in
conformal geometry.  This centres around the notion of conformal blow-up by
means of which a compact Riemannian manifold (or orbifold) with a
marked point is changed into a manifold with a conformally related
$b$-metric.  There is an analogous notion of conformal blow-down for
weakly asymptotically locally euclidean (WALE) spaces (\S\ref{walex}).

We also give the elements of the Fredholm theory of $b$-elliptic
operators. The latter involves us in a short account of the notion of
`polyhomogeneous' functions: such functions arise naturally and
inevitably in the  study of $b$-differential operators 
and should be thought of as an extension of the idea of a function
that is smooth up to the boundary.

Finally in this section, we give a careful description of the process
of gluing a pair of Riemannian $b$-manifolds, and the construction of
a suitable $b$-metric on their join. Thus all the material in this
section is well known to the right people; it is necessary to
summarize it here in the interests of making the present paper
self-contained and fixing notation that will be used throughout.

\subsection{The $b$-category}
\label{bcat}
We begin with the basic definitions, and then provide some 
motivation. The reader can find a detailed account in \cite{tapsit}.

Let $X$ be a compact $n$-manifold with smooth boundary $\p X$; neither
$X$ nor $\p X$ are assumed connected. When working near $\p X$ it is
convenient to fix a function $x\in \ci(X)$ with values in $[0,2]$,
such that $x(p)=0$ if and only if $p\in \p X$ and $dx(p)\not=0$ for
all $p$ with  $x(p) \in [0,1]$. Such a function is often called a {\em
boundary defining function}. Here and subsequently, $f\in C^\infty(X)$
means that $f$ is smooth up  to the boundary of $X$.

Near a boundary point $p$, one can 
introduce adapted coordinate systems $(x,y_1,\cdots,y_{n-1})$, where 
the $y_j$ are local coordinates near $p$ in $\p X$. Using such 
coordinates we can define the $b$-tangent and cotangent bundles 
$\mbox{}^bTX$ and $\mbox{}^bT^*X$. These are smooth bundles of rank 
$n$ over $X$; over the interior $X^o = X\backslash \p X$ they are
equal to $TX^o$, $T^*X^o$ respectively, but at the boundary, near $p$,
$\mbox{}^bTX$ is
spanned by the elements
\begin{equation} \label{bbasis}
x\frac{\p}{\p x}, \frac{\p}{\p y_1},\ldots, \frac{\p}{\p y_{n-1}}
\end{equation} 
and dually $\mbox{}^bT^*X$ is spanned by the elements
$$
\frac{dx}{x}, dy_1,\ldots, dy_{n-1}.
$$
It is easy to see that $\mbox{}^bTX$ and $\mbox{}^bT^*X$ are smooth up 
to the boundary of $X$; the basic idea of the $b$-category is to use 
these bundles in place of $TX$ and $T^*X$ in the development of 
differential analysis and geometry for manifolds with boundary. 

A basic example is the definition of a Riemannian $b$-metric. By 
definition, this is just a positive-definite inner product
on $\mbox{}^bTX$, smooth up to the boundary.   We shall not need the 
most general such metric, but only `exact' $b$-metrics, which take the 
form
\begin{equation} \label{exact}
\mbox{}^bg=    \frac{dx^2}{x^2} + h(x,y),
\end{equation}
where $h(x,y)$ is a symmetric tensor such that $\mbox{}^b g$ is
everywhere positive-definite.  To begin with we assume that $h(x,y)$
is $\ci$ up to $\p X$, but we shall soon need to allow $h$ to be
merely polyhomogeneous (see \S\ref{sphg}). We now give some examples to
explain how $b$-metrics arise naturally in conformal geometry. 

\subsubsection{Example: conformal blow-up} \label{confb}
Let $(\bX,\bg)$  be a compact $n$-dimensional Riemannian manifold without
boundary, let $0\in \bX$ be any point. In geodesic polar coordinates
centred at $0$, we have
$$
\bg = dr^2 + r^2(d\omega^2 + r^2\eta(r,\omega))
$$
where $r$ is geodesic distance from $0$,
$d\omega^2$ is the standard round metric on $S^{n-1}$ and 
$\eta(r,\omega)$ is a family of symmetric tensors on $S^{n-1}$,
uniformly bounded as $r\to 0$.  By multiplying
$g$ by a constant, we can assume that these coordinates are defined
for $0 < r < 2$, say.
The (oriented)
blow-up $X$ of $\bX$ at $0$ is a smooth manifold with boundary $X$
canonically identified with the unit sphere-bundle of $0$ in
$\bX$. The pull-back  $x$ of $r$ is then a boundary defining function and
the pull-back of $g = r^{-2}\bg$ is the $b$-metric
$$
g = \frac{dx^2}{x^2} + d\omega^2 + x^2\eta
$$
where $\eta$ is smooth up to $\p X$.  We call $(X,g)$ the {\em
conformal blow-up} of $\bX$ at $0$.

\subsubsection{Generalization: conformal blow-up of orbifolds}

Recall \cite{baily} that an $n$-dimensional 
orbifold $\bX$ is defined analogously to a manifold,
but a neighbourhood $U_p$ of $p\in \bX$ is homeomorphic to $\RE^n/a_p$,
where $a_p:\Ga_p \times \RE^n \to \RE^n$ is an effective action of the
finite group $\Gamma_p$ (the `local isotropy group') on
$\RE^n$. The covering map $\RE^n \to U_p$ is called a {\em local
uniformizing chart} centred at $p$.

If $\Ga_p= \{1\}$ then $p$ is a {\em smooth point} of
$\bX$, otherwise $p$ is a {\em  singular} point. The set of all
singular points of $\bX$ is denoted 
$\bX_{\rm sing}$. 

{\em In this paper we shall use the term `orbifold' to mean`orbifold
with isolated singular points'.}  Then for each point $p$,
$a_p(\gamma)$ only fixes $0$ for each $\gamma\not=1$ in $\Ga_p$. 

A Riemannian metric $\bg$ on $\bX \backslash \bX_{\rm sing}$ is called
a {\em smooth orbifold metric}, and $(\bX,\bg)$ is called a Riemannian
orbifold, if the pull-back of $\bg$ to a local uniformizing chart
extends smoothly to $0\in \RE^n$ and $a_p$ acts by isometries. In
particular $a_p$ gives a representation of $\Ga_p$ in $O_n$, the
orthogonal group of $T_0\RE^n$. 

By working in a local uniformizing chart one sees that the conformal
blow-up $(X,g)$ of a singular point $p$ of $\bX$ can be defined as in
\S\ref{confb}. The only difference is that  $\p X$ is canonically identified
with the spherical space-form $S^{n-1}/a_p$.

\subsubsection{Example: conformal blow-down of asymptotically 
euclidean spaces} \label{walex}
Let $(M,g)$ be a Riemannian manifold that is
{\em weakly 
asymptotically locally Euclidean} (WALE). By this we mean that there exists a 
compact subset $K\subset M$ and a diffeomorphism $\phi:(\RE^n 
\backslash B)/a_\infty \to M\backslash K$ (where $B$ is some closed ball 
in $\RE^n$ and $a_\infty:\Ga\times \RE^n \to \RE^n$ is an action of the
finite group $\Ga$) such that
$$
|g - g_0| = O(r^{-k}),\;\;|\nabla_0^j g| = O(r^{-j-k})
$$
where $g_0$ is the Euclidean metric, $r$ is the distance from the 
origin of $\RE^n$ and $\nabla_0$ is the covariant derivative of the 
Euclidean metric.  This definition is mainly of interest when $k\geq 
2$. The term `asymptotically locally euclidean' (ALE) has become
standard for the very strong decay with
$k=4$.  Writing the metric in polar coordinates again, we have
$$
g = dr^2 + r^2 (d\omega^2 + r^{-k}\eta)
$$
where $|\nabla^j_0 \eta| = O(r^{-j-k})$ for $j=0,1,2,\ldots$ and 
$\eta$ is again a family of symmetric 2-tensors on $S^{n-1}/a_\infty$. 
Set $x = r^{-1}$ and we obtain
$$
x^2 g = \frac{dx^2}{x^2} + d\omega^2 + x^k\eta,
$$
a $b$-metric on the radial compactification of $M$. We refer to this 
process as conformal blow-down of a WALE space.

\subsubsection{Remark} Observe that from this `$b$' point of view, the conformal 
blow-up of a point and the conformal blow-down of the infinity of a 
WALE space both look exactly the same; a manifold with boundary a 
spherical space-form, equipped with an exact $b$-metric which is the 
standard metric of constant curvature at the boundary.

\subsubsection{Remark} 
By the change of variables $t= \pm\log x$ , $dt = \pm dx/x$, 
the interior of a $b$-manifold with $b$-metric becomes a manifold with 
a cylindrical  end diffeomorphic to $\p X\times (t_0,\infty)$ or 
$\p X\times (-\infty,t_0)$ ($t_0$ some constant), depending on the sign 
chosen. In these coordinates, a $b$-metric becomes a metric which 
approaches a Riemannian product metric on $\p X \times \RE$
at an exponential rate in $t$, 
with similar estimates on all derivatives of the metric.   We shall 
use this change of variables to describe the process of gluing 
arbitrary $b$-manifolds in \S\ref{glueb}. We shall also often confuse
a manifold with cylindrical ends with a $b$-manifold (even though this
runs counter to the idea of the $b$ category, which is to replace
non-compactness by degeneracy at the boundary!).

In this section we have tried to show that there is a precise sense in
which the $b$-category unifies orbifolds and WALE spaces.  Equally
important is that there is 
a good Fredholm theory for partial differential operators naturally
associated to $b$-metrics. This will now be outlined.

\subsection{On $b$-differential operators}
\label{fully}

Throughout this section, $X$ is a $\ci$ manifold with boundary $\p X$,
$x\geq 0$ is a boundary defining function, $t= -\log x$, and $X^o = X
\backslash \p X$.

\subsubsection{Definition}  A $\ci$ $b$-differential operator
$P:\ci(X,E) \to \ci(X,F)$, where $E$ and $F$ are two vector bundles
over $X$, is a differential operator which can be written locally in
the form $P = p(x,y;x\p_x,\p_y)$ where $p$ is smooth in $(x,y)$ and
polynomial in $x\p_x$ and $\p_y$ (cf.\ \eqref{bbasis}).

 Any operator
`naturally associated to' a $b$-metric---for example the Laplacian or
Dirac operator---will automatically be a $b$-differential
operator. Fredholm theory for such operators has been developed by Lockhart
and McOwen 
\cite{lockhart-mcowen}, and in much more detail by Melrose and 
Mendoza \cite{tapsit} and is the main technical tool needed to prove
our  gluing theorems.  

Given a $b$-differential operator $P$, there is a canonically
associated {\em indicial operator} $I(P)$ which is given locally by
$I(P) = p(0,y;\p_t,\p_y)$, regarded as a $t$-invariant differential
operator on the cylinder $\p X \times \RE$.  In terms of the
cylindrical model of $X^o$, all coefficients of $P -I(P)$ decay
exponentially as $t \to \infty$. 

\subsubsection{Definition} Let $P:\ci(X,E) \to \ci(X,F)$ be a
$b$-differential operator of order $m$. Then
$$\spb(P)  =
\{\lambda \in \CX: \mbox{ there exists }u(y)\not=0\mbox{ such that }
I(P)(e^{i\lambda t}u(y))=0\}.
$$
In other words $\spb(P)$ is the set of complex numbers for $\lambda$ for
which $I(P)$ has a non-trivial exponential solution with exponent
$i\lambda$. Notice that here $E$ and $F$ have implicitly been
trivialized in the $t$ direction along the infinite cylinder $\p X\times
\RE$: to be precise, $E$ and $F$ have been identified with bundles
pulled back from $\p X$ by the projection $\p X\times \RE \to \p X$.

\subsubsection{}
\label{disctsp} It is a basic fact that if $P$ is elliptic over
$X^o$, then $\spb(P)$ is a discrete set which meets every horizontal
strip $\{a < \Im(\lambda)  < b\}$ in a finite
number of points.  For the translation-invariant operator $I(P)$,
real elements of $\spb$ are obstructions to its invertibility in $L^p$
spaces:

\begin{propn} The elliptic operator $I(P)$ extends to a bounded map
$$
L^p_{k}(\p X\times \RE,E) \to L^p_{k-m}(\p X\times \RE,F)
$$
for each $p$ and $k$. This map is invertible if and only if
$\spb I(P)\cap \RE = \emptyset$.
\label{tinbt}\end{propn}

Here and below $L^p_k(X)$ is the Sobolev space of functions $u$ 
such that $Pu \in L^p$ for every $b$-differential operator of order
$\leq k$ and the measure used to define $L^p$ is 
that defined by any $\ci$ $b$-metric on $X$. For the cylinder $\p
X\times \RE$ an example of such a metric would be a ($t$-invariant) Riemannian product
metric.

\subsubsection{Fully elliptic $b$-differential operators} The
$b$-differential operator $P$ is said to be {\em fully elliptic} if
$P$ is elliptic in $X^o$ and $\spb(P)\cap \RE = \emptyset$. In view of
the proposition, the latter is saying that $P$ is `invertible at
$\infty$ (or at $\p X$)'. It turns out that an elliptic
$b$-differential operator is Fredholm in $L^p$ if and only if it is
fully elliptic.  Before stating the theorem which summarizes the
mapping properties of fully elliptic $b$-differential operators, we
need to introduce polyhomogenous functions.

\subsubsection{Polyhomogeneity} \label{sphg} If $X$ is the conformal blow-up at a
point of $\bX$, then the pull-back $u$ to $X$ of a $\ci$ function on $\bX$
is smooth and in particular has an asymptotic expansion $u \sim \sum_{j = 0}^\infty
u_j(y)x^j$ near the boundary, where $u_j\in \ci(\p X)$ for each $j$.

A polyhomogenous (phg) function is $\ci$ in $X^o$ and has a
similar asymptotic expansion near $\p X$, but more general powers
of $x$ can appear, as well as polynomials in $\log x$. The set of
powers that can occur is called an {\em index set} $\CI$ and must be a
discrete subset of $\CX$ having the additional property that if $z_j
\in \CI$ and $|z_j| \to \infty$, then $\Re(z_j) \to +\infty$.  Given
an index set $\CI$,
$u\in \ci(X^o)$ is called polyhomogeneous (with respect to $\CI$) if 
\begin{equation} \label{phgeq}
u \sim \sum_{z\in \CI} x^zu_z(y,\log x)
\end{equation}
where for each $z$, $u_z(y,t)$ is $\ci$ in $y$ and  polynomial in
$t$. The symbol $\sim$ in \eqref{phgeq} is meant in the following
strong sense: if
$$
u_N = \sum_{z\in \CI,\Re(z) \leq N} x^zu_z(y,\log x),
$$
then $u-u_N\in C^N(X)$ and all derivatives of order $\leq N$ of
$u-u_N$ vanish at $\p X$. For more details, see
\cite[\S5.10]{tapsit}; there, however, a more refined notion of index
set is used, designed to keep track of the degrees of the polynomials
$u_z(y,\cdot)$. We have chosen to ignore this refinement here.

We are now ready to summarize the essential properties of fully
elliptic operators on $X$.
\begin{thm} \label{poly}
Let $P:\ci(X,E) \to \ci(X,F)$ be a fully elliptic
$b$-differential operator of order $m$. Then
\begin{enumerate}
\item[(i)] $P$ extends to a Fredholm map $L^p_{k}(X,E) \to L^p_{k-m}(X,F)$
      for every $p$ and $k$, with index independent of $p$ and $k$.

\item[(ii)] If $u\in \ker(P)$ then $u$ is $\ci$ in $X^o$ and is phg
relative to the index set $\CI = i\spb(P)\cap \{\Re z > 0\}$.

\item[(iii)] The cokernel of $P$ can be identified with the $L^p$ kernel of
      $P^*$, where $P^*$ is the $L^2$-adjoint of $P$ with respect to
      a $b$-metric on $X$. In particular $PL^p_k(X,E)$ can be
      complemented by phg sections of $F$.
\label{frpk}
\end{enumerate}
\end{thm}

For a proof see \cite{tapsit} or \cite{lockhart-mcowen}. 

\subsubsection{Conjugation and weights}
\label{cjw} Let $\delta$ be a real number
and let 
$P^{(\delta)}= x^{\delta} P x^{-\delta} = e^{-\delta t}Pe^{\delta t}$. 
Then 
$\spb(P^{(\delta)}) = \spb(P) - i\delta$ and so if $P$ is elliptic,
then $P^{(\delta)}$ will be fully elliptic for all $\delta \not\in
\Im\spb(P)$. By the remark at the beginning of \S\ref{disctsp}, 
{$P^{(\delta)}$ is fully elliptic for all
but a discrete set of real values $\delta$.  The index of
$P^{(\delta)}$ is locally constant in $\delta$ and jumps as $\delta$
passes through a point in $\Im\spb(P)$.

An equivalent formulation of this observation is that an elliptic
$b$-operator $P$ defines a Fredholm operator between {\em weighted
Sobolev spaces} 
$$
e^{\delta t}L^p_{k}(X,E) \to e^{\delta t}L^p_{k-m}(X,F)
$$
for all $\delta \not\in \Im\spb(P)$.  Note carefully, however, that
the weighted Fredholm alternative (analogue of Theorem~\ref{frpk}
(iii)) identifies the cokernel of this map with the kernel of the map
$$
e^{-\delta t}L^p_{k}(X,F) \to e^{-\delta t}L^p_{k-m}(X,E).
$$

We shall need to make use of
these ideas, for $0\in \spb$ for the linearization of the conformal ASD
equations.

\subsubsection{Remark} If $\p X= Y_1\cup\ldots\cup Y_n$ where the $Y_j$
are connected
then $\spb= \spb^{(1)}\cup \ldots \cup \spb^{(n)}$ where $\spb^{(j)}$
is the contribution from $Y_j$. Then one can shift these pieces
independently of each other by conjugating by a function which is
equal to $e^{\delta_jt_j}$ near $Y_j$. We shall not develop a
systematic notation for this situation.

\subsubsection{Polyhomogeneous $b$-differential operators}

Since polyhomogeneous functions are at least as natural on a manifold
with boundary as functions that are $\ci$ up to  the boundary, it is
natural to widen the class of operators we consider by allowing their
coefficients to be polyhomogeneous relative to some index set
${\mathcal J}$ and continuous up to the boundary.
Such operators arise naturally as operators canonically associated to
a {\em phg $b$-metric} i.e. a metric $\mbox{}^bg$ of the form
\eqref{exact} where $h(x,y)$ is phg and continuous up to the
boundary. The results of this section, in particular
Theorem~\ref{frpk}, go through in this case,  the only difference
being in part (ii) where $u$ will now have a phg expansion relative to
the index set $\CI \cup  {\mathcal J}$, $\CI$ being as before.

\subsection{On gluing $b$-manifolds}
\label{glueb}
We shall now explain how to glue $b$-manifolds across (a part of)
their boundaries. The construction is complicated slightly by the need
to keep track of orientations.  The main points of the discussion are
contained in \S\ref{xrdef}, \S\ref{gbr} and \S\ref{firstp}.

\subsubsection{Data}
\label{data} For $j=1, 2$ let $X_j$ be a smooth, oriented, compact
$n$-manifold with boundary and suppose that an oriented boundaryless manifold
$Y$ occurs in $\p X_j$ with each orientation:
$$ Y \subset \p X_1,\;\; -Y \subset \p X_2. $$
(In what follows, we shall not indicate the orientation of $Y$ unless
orientation issues are being discussed.) 
Let $x_j$ be defining functions for $Y \subset \p X_j$ and assume that
$x_j=1$ near $\p X_j \backslash Y$. Let $t_1 = -\log x_1$, $t_2 = \log
x_2$; then there exist open sets $U_j \subset X_j^o$ diffeomorphic to
cylinders
$$
U_1 = Y\times \{0 < t_1 < \infty\},\;\;U_1 = Y\times \{-\infty < t_2 < 0\},
$$
such that a sequence of points tending towards $Y$ corresponds to
$|t_j| \to \infty$.  It is important to notice that $U_j$ 
inherits from $X_j$ is $dt_j\wedge \mbox{or}_Y$ where
$\mbox{or}_Y$ is a given orientation of $Y$.

\subsubsection{Definition}
\label{xrdef} The manifold $X_\rho$ is defined by
truncating the $U_j$ at $\pm t_j = \rho$ and identifying the
boundaries $Y\times \{t_1=\rho\}$ and $Y\times \{t_2=-\rho\}$.

After the remarks of the previous paragraph, it is clear that $X_\rho$
is oriented, this orientation agreeing with the given orientations on
the $X_j$. Furthermore $X_\rho$ contains a {\em neck} of length
$2\rho$ given by $\{0 < t_1 \leq\rho\}\cup \{-\rho \leq t_2 \leq 0\}$.

\subsubsection{Notation} On $X_\rho$ the function $t$ is defined so
that $t = t_1-\rho$ for $0 \leq t_1 \leq \rho$ and $t= t_2+\rho$ for
$-\rho \leq t_2 \leq 0$, and extended smoothly to $-\rho$ on the rest
of $X_1$ and to equal $\rho$ on the rest of $X_2$. The central slice
of the neck then corresponds  to $t=0$. See Figure~\ref{giraffe}.
\begin{figure}
 \begin{center}
 \begin{tabular}{c}
 {\epsfig{figure=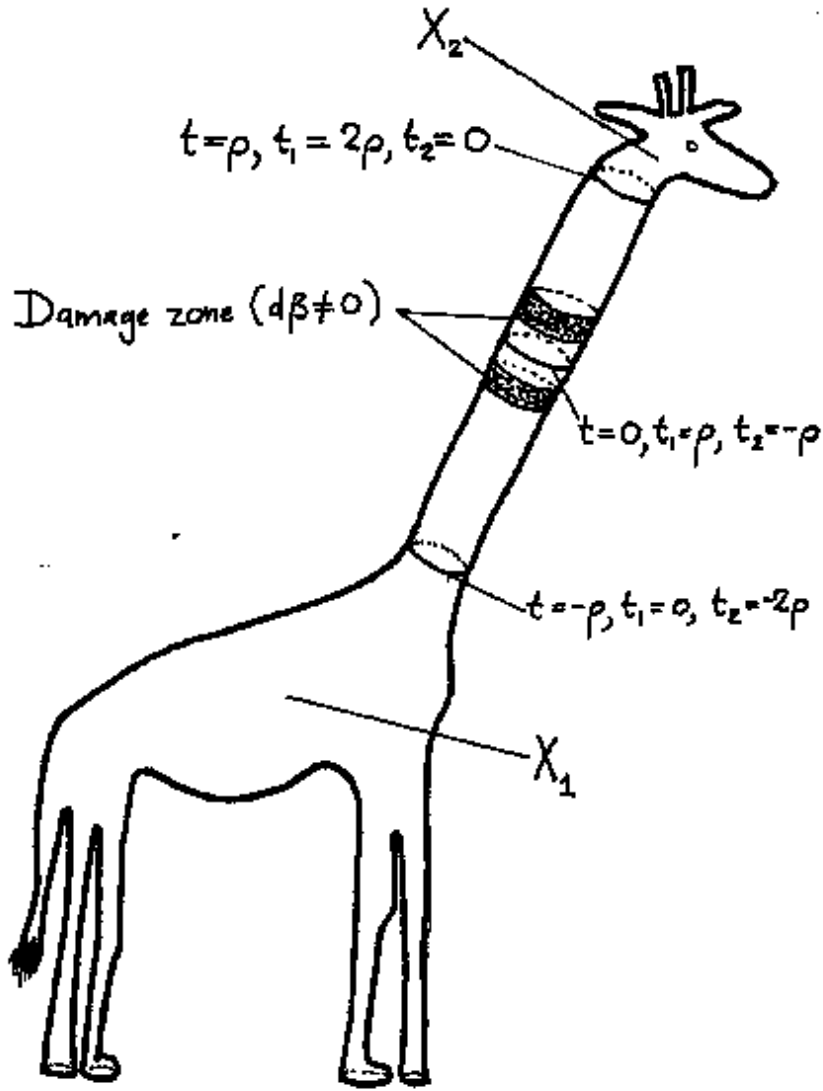,width=14.8cm}}
 \end{tabular}
\caption{The manifold $X_\rho$ with a long neck constructed by gluing the
$b$-manifolds $X_1$ and $X_2$. (Adapted from a drawing by F.\ N.\ Singer.)}
 \label{giraffe}
 \end{center}
 \end{figure}

\subsubsection{Example} Suppose that $X_1$ is the conformal blow-up at a point
$0$, say, of a compact orbifold $\bX$. In particular a neighbourhood
of $0$ in $\bX$ is homeomorphic to the quotient $\RE^4/a$, where
$a:\Ga \times \RE^4 \to \RE^4$ is an action of the finite group $\Ga$
on $\RE^4$.  Suppose further that $X_2$ is the conformal
blow-down of the asymptotic region of a WALE space $\hX$, where a
neighbourhood of $\infty$ is homeomorphic to a neighbourhood of
$\infty$ in $\RE^4/a$.  Then $\p X_1 = Y = S^3/a$ and $\p X_2 = - Y$
and the construction of \S\ref{data} can be applied. The result is a
generalized blow-up of $0\in \bX$, where a small ball centred at $0$
is replaced by the complement of a neighbourhood of infinity in $\hX$.
\label{wale}

\subsubsection{Gluing Riemannian $b$-metrics} 
\label{gbr}
Let $X$ be a $b$-manifold. A $\ci$ Riemannian metric on $X^o$ is
called a polyhomogeneous $b$-metric on $X$ if and only if it has the
form \eqref{exact} where $h(x,y)$ is continuous up to $\p X$ and has a
polyhomogeneous expansion relative to some index set.  Assume now that
the $X_j$ in \S\ref{data} are equipped with such phg $b$-metrics $g_j$
and that near $Y$, $h_1(0,y) =
h_2(0,y)= :h(y)$, say. In terms of the cylindrical parameters $t_j$, the $g_j$
both approach the metric
$$
g_0 = dt^2 + h(y)\mbox{ on }X_0:=\RE\times Y
$$
as $|t_j| \to \infty$. 

We construct a Riemannian metric $g_\rho$ on $X_\rho$ by picking once
and for all a standard non-increasing cut-off function $\beta:\RE \to [0,1]$, such 
that $\beta(t) =1$ for $t\leq -1/2$ but $\beta(t) =0$ for $t\geq 1/2$
and using it to define new metrics 
$$
\wg_{1,\rho} = \beta(t_1 -\rho+1)g_1 + (1 - \beta(t_1 -\rho+1))g_0
$$
on $X_1$ and
$$
\wg_{2,\rho} = \beta(\rho-1- t_2)g_2 + (1 -\beta(\rho-1-t_2))g_0
$$
on $X_2$. These formulae cut off the exponentially decreasing terms in
the asymptotic expansions of the $g_j$ leaving the standard
$t$-independent metric $g_0$ for $|t_j| \geq \rho$. In particular the
identification used in the construction of $X_\rho$ is now an isometry
and we define $g_\rho$ on $X_\rho$ to be equal to $\wg_{1,\rho}$ for
$t\leq 0$ and to be equal to $\wg_{2,\rho}$ for $t\geq 0$.

\subsubsection{First properties of $g_\rho$} \label{firstp} It is clear from the
construction that $g_\rho$ is equal to $g_1$ or $g_2$ for $|t| \geq 2$
and is equal to $g_0$ for $|t| < 1/2$. Even in the {\em damage zone}
$\{-1/2\leq |t-1| \leq 1/2\}$, $g_\rho$ is {\em exponentially close} to
$g_0$.  By this we mean that there exists $\eta >0$ such that
\begin{equation} \label{metdev}
\sup_{K}|g_\rho - g_0| = O(e^{-\eta\rho})\mbox{
as }\rho\to\infty,
\end{equation}
where $K = \{-1/2\leq |t-1| \leq 1/2\}$ and the pointwise norm is that
induced by $g_0$.  In fact, \eqref{metdev} holds for $K = \{|t|\leq
T\}$, for any fixed $T>0$, and similar estimates hold for all
derivatives of $g-g_0$; such estimates follow from the assumed
polyhomogeneous expansions of the $g_j$ near $Y$.

In particular if $P_j$ and $P_\rho$ are differential operators
canonically associated to the metrics $g_j$ and $g_\rho$ then $P_\rho$
is exponentially close to the $P_j$ and hence also to 
$P_0$ on $\{|t| < T\}$ for any
fixed $T$, as $\rho \to \infty$.

\setcounter{equation}{0}
\section{The anti-self-duality equations}

In this section we first review Riemannian and conformal geometry, 
passing in \S\ref{bak4d} to the special case of $4$ dimensions where 
we stay for most of the rest of the paper. In particular, the 
relevant analytical aspects of the ASD equations are given here, the 
main facts being summarized in Proposition~\ref{asdsum}. In \S\ref{hasd}
we study the Hermitian version of the ASD equations from the same 
point of view, listing the main points in Proposition~\ref{hasdsum}. 
Thus experts in $4$-dimensional geometry may well be able to move 
straight to these Propositions, referring back if necessary to the 
earlier parts of this section.

\subsection{Preliminaries on metrics and conformal structures}

\subsubsection{Metrics and curvature}
\label{metcurv}
Let $(X,g)$ be a Riemannian $n$-manifold and let $h$ be 
a $g$-symmetric endomorphism of $TX$. Then the bilinear form
\begin{equation}\label{gh}
    g_h(\xi,\eta):= g(\xi,\eta) + g(h\xi,\eta)\; (\xi,\eta\in 
    C^\infty(X,TX))
\end{equation}
is symmetric and defines a Riemannian metric if $|h| 
:=(\tr(h^2))^{1/2} <1$ at each point of $X$.     Let $\nabla$ and 
$\nabla^h$ denote respectively the metric connections of $g$ and
$g_h$. Then we have
\begin{equation} \label{chcon}
\nabla_h\xi  = \nabla\xi +  Q((1+h)^{-1},\nabla h)\xi
\end{equation}
for any section $\xi$ of $TX$, where $Q(u,v)$ is an
endomorphism-valued $1$-form, bilinear in $(u,v)$.

Hence the curvature $R_h$ of $g_h$ has the form
\begin{equation} \label{chcurv}
R_h = R_0 + d^\nabla Q((1+h)^{-1},\nabla h) + 
Q((1+h)^{-1},\nabla h)\wedge Q((1+h)^{-1},\nabla h).
\end{equation}
where $R_0$ is the curvature of $g$ and $d^\nabla$ is the covariant 
exterior derivative defined by $\nabla$. Separating the linear and
nonlinear terms in \eqref{chcurv}, we obtain
\begin{equation}\label{betcurv}
R_h = R_0 + R'_0[h] + \varepsilon_1(1+h,\nabla h\otimes\nabla h) 
+\varepsilon_2(1+h,h\otimes\nabla\nabla h)
\end{equation}
where $R'_0$ is a linear differential operator and each
$\varepsilon_j$ is real-analytic in the first variable and linear in
the second variable, with coefficients 
depending only upon $g$ and its derivatives.
\subsubsection{Convention}
\label{conv} From now on we shall use
$\varepsilon(u,v)$,  $\varepsilon_1(u,v)$,  $\varepsilon_2(u,v)$,
etc. generically for  `error terms' with the properties just
mentioned. That is, they are real-analytic in the $0$-jet of $u$ near
$u=1$ and linear in the $0$-jet of $v$. 
For example, the $\varepsilon_j$ in
\eqref{betweyl} below are not identical to those in
\eqref{betcurv}. The point is that in order to estimate these
non-linear terms later on, all we shall need to know is their qualitative
dependence upon $h$ and its derivatives.

\subsubsection{Conformal structures}
\label{coco}
We adopt a `modern' approach to conformal structures for which we
claim no originality. For more details of this approach,  the reader
could consult, for example, \cite{dmjc}.

Let $X$ be a
$C^\infty$ manifold and let $\Om = |\Lambda^n T^* X|$
be the bundle of densities on $X$ (\cite[p.\ 148]{h1}). 
Since $\Om$ is an $\RE_+$-bundle, $\Om^w$ has a canonical
meaning for any real $w$.  By a {\em Riemannian conformal structure}
on $X$ we shall mean a suitably normalized positive-definite $\ci$
section of the bundle $S^2 T^*X \otimes \Om^{-2/n}$.  Since the top
exterior power of $\Om^{1/n}TX$ is canonically trivial a possible
normalization is the condition $\det c =1$, but we shall make a
different choice in \S\ref{paramc}. (Here we have written $\Om^{1/n}TX$ for
$\Om^{1/n}\otimes TX$. We continue to omit such tensor product signs below.)
Note that $c$ can also be viewed as a
normalized metric on the {\em weightless tangent bundle}
$\Om^{1/n}TX$.

A positive trivialization or {\em choice of length-scale} 
$\mu$ of $\Om^{-1/n}$ determines a {\em
compatible Riemannian metric} $g_\mu = \mu^{-2}c$ on $TX$; any two
such are {\em conformally related} in the sense that
\begin{equation}\label{metch}
g_{\mu'} = (\mu'/\mu)^2 g_\mu
\end{equation}
where $\mu'/\mu$ is a positive $\ci$ function, so that any two
compatible metrics are related by {\em conformal rescaling}. In
particular, the present approach is equivalent to the more traditional
one in which a conformal structure is taken as a conformal equivalence
class of Riemannian metrics. We shall occasionally 
write $g\in c$ to mean that $g$ is
compatible with $c$, i.e. that $g$ arises from $c$ by a choice of
scale $\mu$ or say that $g$ belongs to the conformal class of $c$.

The length-scale $\mu$ also
defines a metric $|\cdot|_\mu$ on the tensor bundle
$\Om^{-w/n}TX^{\otimes j}\otimes T^*X^{\otimes k}$
 (and any sub- or quotient bundle) by
the formula
\begin{equation}\label{mtdef}
|s|_\mu : = |\mu^{-w-j+k}s|_c
\end{equation}
(generalizing \eqref{metch}).
This bundle is therefore said to have {\em conformal weight} $w+j -k$,
because of the formula
\begin{equation}\label{wtch}
|s|_{\mu'} = (\mu'/\mu)^{w+j-k}|s|_{\mu}.
\end{equation}
Note $c$ determines an identification $\Om^{1/n}TX \to \Om^{-1/n}T^*X$,
the conformal version of index-lowering. Similarly indices are raised
with $c^{-1}$. Note that these operations {\em preserve conformal
weight}.

\subsubsection{Definition} A bundle $E$ of conformal weight $0$
is called a {\em weightless} bundle. 

If $E$ is a weightless bundle then a conformal structure $c$ defines a
genuine metric on $E$.  It is clear that any bundle associated to the
tangent bundle of $X$ can be written uniquely in the form
$\Om^{-a/n}E$, where $E$ is weightless. We shall frequently write
bundles in this way in situations where it is necessary to keep track
of conformal weights.

\subsubsection{Parameterization of conformal structures}
\label{paramc}
Let  $(X,c)$ be a conformal $n$-manifold with 
$\det c =1$. 
If $h$ is a $c$-symmetric endomorphism of $TX$, we define
\begin{equation}\label{ch}
    c_h(\xi,\eta):= c(\xi,\eta) + c(h\xi,\eta)\; (\xi,\eta\in 
    C^\infty(X,TX)).
\end{equation}
For suitable $h$ (in particular if $|h| < 1$ at each point) then $c_h$
will be positive-definite. 
We choose to normalize $c_h$ by the requirement that $h$ be
trace-free, rather than $\det c_h =1$. To first order in $h$, these
conditions agree.

Denote by $E^1$ the bundle of $c$-symmetric trace-free endomorphisms of
$TX$; then $E^1$ is weightless and using $c$ can be identified
with $\Om^{-2/n}S^2_0T^*X$, where $S^2_0=$ trace-free part of the
symmetric square.

In terms of $E^1$, we can summarize this paragraph as follows:
there is an open 
neighbourhood $B$ of the zero-section in $C^\infty(X,E^1)$, 
containing all $h$ with $\sup_X |h| < 1$, such that $h \mapsto c_h$ is 
a diffeomorphism of $B$ with the set of all conformal structures on $X$.

\subsubsection{Conformal invariance}
\label{sconfinv}
It is well known in conformal geometry that certain differential
operators, which {\em a priori} depend upon a Riemannian metric, are
unchanged by conformal rescaling. If such a {\em conformally
invariant}  operator $P:C^\infty(\Om^{-a/n}E)\to C^\infty(\Om^{-b/n}F)$ is of order
$k$, say, 
where $E$ and $F$ are weightless bundles,  then its symbol provides a
map 
$$
\Om^{-a/n}S^kT^*X\otimes E
= \Om^{-a/n+k/n}S^k(\Om^{-1/n}T^*) \otimes E \to \Om^{-b/n} F
$$
which must be weightless. Accordingly we must have  $b= a -k$.

The formal adjoint of a conformally invariant operator defines a
conformally invariant operator but the conformal weights may
change. Indeed it is a familiar fact in the analysis literature (e.g.\
\cite[p.\ 93]{h3} that a
differential operator $P:C^\infty(\Om^{1/2}E) \to
C^\infty(\Om^{1/2}F)$ has a formal adjoint $P^*:C^\infty(\Om^{1/2}F^*)
\to C^\infty(\Om^{1/2}E^*)$ independent of any metric. If we replace
$\Om^{1/2}E$  and $\Om^{1/2}F$  by $\Om^{-a/n}E$ and $\Om^{k/n-a/n}F$
where $E$ and $F$ are now assumed weightless, then using $c$ to identify
$E$ with $E^*$ and $F$ with $F^*$, we obtain
$$
P^*:C^\infty(\Om^{(1 + (a-k)/n)}F) \to C^\infty(\Om^{(1 + a/n)}E)
$$
which is conformally invariant if $P$ is conformally invariant.

\subsection{Background to $4$-dimensional geometry}
\label{bak4d}
$4$-dimensional Riemannian geometry is enriched by the existence of
the special isomorphism 
${\mathfrak so}(4) = {\mathfrak so}(3)\oplus {\mathfrak so}(3)$. 
 The geometric counterpart of this algebraic fact
is the decomposition 
\begin{equation} \label{2fdec} \Lambda^2=
\Lambda^+\oplus \Lambda^- \end{equation}
 for
$2$-forms on a $4$-dimensional vector space equipped with a metric and
orientation.  As far as the present work is concerned, the main
consequence of this is the presence of the anti-self-duality equations
for curvatures on an oriented $4$-manifold. We start, however, with
some basic algebra.

\subsubsection{Algebraic preliminaries} \label{algpr} Consider euclidean space
$\RE^4$ with its standard metric and orientation $dx_0\wedge dx_1
\wedge dx_2  \wedge dx_3$.  The purpose of this paragraph is to
explain the invariant isomorphism $S^2_0\RE^4= \Lambda^+\otimes
\Lambda^-$, which will often be used below.

Consider the standard bases
$$
e_1 = dx_0\wedge dx_1 + dx_2\wedge dx_3,\;
e_2 = dx_0\wedge dx_2 + dx_3\wedge dx_1,\;
e_3 = dx_0\wedge dx_3 + dx_1\wedge dx_2,
$$
of $\Lambda^+$ and
$$
\be_1 = dx_0\wedge dx_1 - dx_2\wedge dx_3,\;
\be_2 = dx_0\wedge dx_2 - dx_3\wedge dx_1,\;
\be_3 = dx_0\wedge dx_3 - dx_1\wedge dx_2,
$$
of $\Lambda^-$. Using the metric, $e_r$ and $\be_r$
operate on $\RE^4$ as orthogonal complex structures $I_r$ and $\bI_r$,
respectively, corresponding to the left and right action of the
quaternions $i$, $j$ and $k$. In particular 
$$I_1^2 = I_2^2 = I_3^2 = I_1I_2I_3= -1 $$
but
$$I_r\bI_s = \bI_sI_r \mbox{  for all }r\mbox{ and }s.
$$

Now let $A$ be an endomorphism of $\RE^4$. The irreducible components
of $A$ consist of
\begin{itemize}
\item $\tr(A)\in \RE$,
\item $\sum \tr(I_rA)e_r \in \Lambda^+$,
\item $\sum \tr(\bI_rA)\be_r \in \Lambda^-$,
\item $\sum \tr(I_r\bI_s A)e_r\otimes \be_s \in \Lambda^+\otimes\Lambda^-$,
\end{itemize}
and the maps implied by these formulae are equivariant. It is easy to
see, using standard properties of the trace, that the 
last of these is indeed an isomorphism of the space of symmetric
trace-free endomorphisms of $\RE^4$ with $\Lambda^+\otimes \Lambda^-$.

{\em For the rest of this section, $(X,c)$ will be an oriented
conformal $4$-manifold, and $g$ will be a Riemannian metric in the
conformal class $c$.}

\subsubsection{The anti-self-duality equation for conformal
structures} The curvature $R$ of $g$, viewed as a symmetric
\label{asd-basics}
endomorphism of $\Lambda^2$  decomposes as
\begin{equation}\label{curv1}
    R = \begin{pmatrix}
    W^+ +s/12 & \Phi\cr
    \Phi^t & W^- + s/12\cr
    \end{pmatrix} : 
\begin{pmatrix} \Lambda^+ \cr \Lambda^- \cr \end{pmatrix}
    \to
    \begin{pmatrix} \Lambda^+ \cr \Lambda^- \cr \end{pmatrix}
    \end{equation}
where $\Phi$ is
the trace-free
part of the Ricci tensor of~$g$, viewed as a section of
$\Lambda^-\otimes\Lambda^+$ as in \S\ref{algpr}, $s$ is the scalar curvature and 
$W^+$ and $W^-$ are respectively the self-dual and anti-self-dual 
parts of the {\em Weyl curvature} $W$.  The metric $g$ is said to be 
{\em conformally ASD} if $W^+(g)=0$. Since $W= W^++W^-$ and \eqref{2fdec}
are  conformally invariant, we also write $W^+(c)=0$ and call $c$ an
ASD conformal structure.

\subsubsection{The  deformation complex}
\label{ASDdef}  If $W^+(c) =0$ then to first order in $h$ the
condition $W^+(c_h) =0$ is equivalent to a conformally invariant
differential equation which will be denoted $D_c h =0$.  This operator
is part of the {\em deformation complex} for the conformal ASD equations,
	\begin{equation}
	\label{deform-cx}
C^\infty(X,\Om^{-1/4}E^0) \xrightarrow{L_c} C^\infty(X,E^1) 
\xrightarrow{D_c} C^\infty(X,\Om^{1/2}E^2).
	\end{equation}
of conformally invariant operators,  where
\begin{equation} \label{wtedbundles}
\Om^{-1/4}E^0 = TX,\; E^1 = \Om^{-1/2}S^2_0 T^*X,\; E^2 = \Om^{-1}
S^2_0\Lambda^+T^*X.
\end{equation}
In particular each of the $E^j$ is weightless.  As we have said, $D_c$ is the
linearization of the map $h \mapsto W^+(c_h)$ and is a second-order
operator; $L_c$ gives the action of
infinitesimal diffeomorphisms on conformal structures, so that
$L_c\xi$ is equal to the trace-free part of the Lie derivative of $c$
along $\xi$.

The operator $D_c$ is easier to understand after use of the isomorphism
 $S^2_0(\Om^{-1/4}T^*X) = \Lambda^-(\Om^{-1/4}T^*X)\otimes
\Lambda^+(\Om^{-1/4}T^*X)$ described in \S\ref{algpr}.
Suppressing powers of $\Om$ for the moment,
recall the existence of a natural second-order elliptic operator
$d^+d^*:\ci(\Lambda^-) \to \ci(\Lambda^+)$.  This can be coupled to
any vector bundle $V$ with connection $A$ to yield an operator
$d^+_Ad^*_A:\ci(\Lambda^-\otimes V) \to \ci(\Lambda^+\otimes
V)$. Then the second-order part of $D_c$ is the composite of this
map (with $V = \Lambda^+$) and projection $\Lambda^+\otimes\Lambda^+
\to S^2_0\Lambda^+$.  $D_c$ also has a zeroth-order term given by
multiplication by the trace-free part of the Ricci tensor.

The composite $D_cL_c$ vanishes if $c$ is ASD 
and then the {\em deformation cohomology groups} $H^*_c$ are
defined and have a standard role \cite{king-kotschick} in describing the local
properties of the moduli space of ASD conformal structures near $c$.
\eqref{deform-cx} is an elliptic complex, so the $H^*_c(X)$ are
finite-dimensional vector spaces if $X$ is a compact manifold (or
orbifold) without boundary.

\subsubsection{The full ASD equations}
The self-dual Weyl tensor $W^+_h:= W^+(c_h)$ is obtained by orthogonal
projection of $R_h$ onto 
$\Om^{1/2}E^2_h = \Om^{-1/2} S^2_0\Lambda^+_h$, where
$\Lambda^+_h$ is the bundle of self-dual $2$-forms for
$c_h$.  Since this orthogonal projection depends real-analytically
upon the $0$-jet of $(1+h)$, 
we have a formula analogous to \eqref{betcurv}:
\begin{equation} \label{betweyl}
    W^+_h = W^+_0 + D_c h + \varepsilon_1(1+h,h\otimes\nabla\nabla h) +
    \varepsilon_2(1+h,\nabla h\otimes \nabla h) + 
    \varepsilon_3(1+h,h\otimes h)
\end{equation}
(Recall the convention regarding the $\varepsilon_j$ explained in \S\ref{conv}.)

It is convenient to remove the dependence of the 
target space upon $h$ by projecting orthogonally from $E^2_h$ onto 
$E^2$ (recall that both $E^2$ and $E^2_h$ are subbundles of 
$\Om^{-1}S^2\Lambda^2$).
In this way we obtain a nonlinear map
$F(h): B \to C^\infty(X,\Om^{1/2}E^2)$, 
such that $F(h)=0$ if and only if $c_h$ 
is conformally ASD.  Summing up,
\begin{propn}\label{asdsum} Let $(X,c)$ be an oriented conformal
$4$-manifold and let the $E^j$ be as in \eqref{wtedbundles}. Then
 there exists a $\ci$ map  $F: B \to C^\infty(X,\Om^{1/2}E^2)$, 
where $B:= \{h\in \ci(X,E^1): \sup_X|h| < 1\}$, with an expansion of
the form
\begin{equation}\label{Fdef}
    F(h) = W^+_0 + D_c h + \varepsilon_1(1+h,h\otimes\nabla\nabla h) +
    \varepsilon_2(1+h,\nabla h\otimes \nabla h) +
    \varepsilon_3(1+h,h\otimes h)
\end{equation}
such that $F(h) =0$ if and only if $c_h$ is ASD. In \eqref{Fdef},
$W^+_0$ is the self-dual Weyl tensor of $c$, $D_c$ is the
linearization of $h \mapsto W^+(c_h)$, and the $\ve_j$ are nonlinear
terms which conform to Convention~\ref{conv}.
\end{propn}
\label{sFdef}

\subsection{Hermitian-ASD conformal structures}
\label{hasd}
Let $(X,J)$ be a complex surface and let $c$ be a $J$-hermitian
conformal structure. We shall repeat our discussion of deformations of
ASD conformal structures now preserving the $J$-hermitian condition.

\subsubsection{Hermitian deformations of $c$} In the parameterization
of deformations of $c$ by endomorphisms $h$ (\S\ref{paramc}) it is
easily checked that $c_h$ is $J$-hermitian if and only if $[h,J]
=0$. Working equivalently with $\Lambda^-\otimes\Lambda^+$ (and
forgetting about conformal weights for now) the $J$-hermitian
deformations correspond to elements of the form $u\otimes\omega$ where
$u\in \Lambda^-$ and $\omega$ is the fundamental $2$-form 
defined by $c$ and $J$. 
\label{paramcj}
\subsubsection{The conformal weight of $J$} For reasons that will emerge
in a moment, it is convenient to regard $\omega$ not as
weightless but as a section of the bundle
$\Om^{-1/4}E^2_J:=\Om^{-3/4}\Lambda^+$.  Then $J$ defines not an
endomorphism of $TX$ but a map $\Om^{-w/4}TX \to \Om^{-(w+1)/4}TX$. (In
the presence of $c$, $\omega$ and $J$ are interchangeable.)  Given $c$
and $\omega\in \Lambda^+$ such that $\omega$ is nowhere $0$ in $X$,
there is a unique choice of gauge such that $J_\mu:=\mu^{-1}\omega$
satisfies $J_\mu^2 = -1$, so there is no loss in giving $\omega$ this
strange conformal weight. Having done so,  we are forced to take 
$u\in \ci(\Om^{1/4}E^1_J := \Om^{-1/4}\Lambda^-)$ so that 
$h=u\otimes \omega$ is weightless.  We now explain why this choice of
conformal weight is advantageous.
\label{crfwt}

\subsubsection{Scalar-flat K\"ahler metrics}
\label{sfk}
It is well known \cite{boyer,vai} that a 
K\"ahler metric $g$ on a complex surface
$(X,J)$  is conformally ASD if and
only if its scalar curvature is $0$; we refer to such $g$ as {\em
scalar-flat K\"ahler}. For any K\"ahler metric $g$, we have $\nabla
J=0$; the complex structure is parallel for the metric connection. An
alternative characterization is that the associated K\"ahler form
$\omega$ should satisfy the {\em twistor equation}
\begin{equation} \label{te}
T\omega = 0,\;\; T\omega = (\nabla\omega)_0
\end{equation}
where the subscript $0$ denotes trace-free part, that is, the image
by the projection
$$
\Lambda^1\otimes\Lambda^+ \to (\Lambda^1\otimes\Lambda^+)_0
$$
(whose kernel is $\Lambda^1$).  Indeed, a short calculation shows that
if $T\omega =0$ and $|\omega| = 1$, then $\nabla\omega = 0$.

On the other hand, $T$ is a  conformally invariant operator
$$
T: \ci(\Om^{-1/4}E^2_J) \to \ci(( \Om^{-3/4}\Lambda^1\otimes \Lambda^+)_0)
$$
so that we can regard a pair $(c,\omega)$ where $c$ is a conformal
structure and $\omega$ is a nowhere-vanishing solution of the twistor
equation \eqref{te},
as being a {\em conformally invariant description}
of a K\"ahler metric on a $4$-manifold.  Furthermore, $c$ is ASD if
and only if the K\"ahler metric is scalar-flat K\"ahler.

As in \S\ref{crfwt}, the K\"ahler representative of $c$ is
fixed by the unique scale which gives $\omega$ constant length equal
to $1$ at each point of $X$.

Finally we note that any Hermitian-ASD conformal structure on a
compact complex surface $X$  of K\"ahler type (equivalently $b_1(X)$
even) is conformally equivalent to a scalar-flat K\"ahler metric 
\cite{boyer,vai}. Thus looking for Hermitian-ASD conformal structures
is, under suitable topological conditions, {\em equivalent} to looking
for scalar-flat K\"ahler metrics. We shall apply this observation to construct
scalar-flat K\"ahler metrics on multiple blow-ups of $\CX^2$ in 
Theorem~\ref{resF}.

\subsubsection{Deformation theory of Hermitian-ASD conformal
structures}

Since hermitian deformations are parameterized by a bundle of rank
$3$, it might be expected that the ASD equations are now
overdetermined. This is not the case: it is shown in \cite{boyer} that
if $c$ is Hermitian, then $W^+(c)$ lies in the image of the map 
$$
\Om^{3/4}E^2_J \stackrel{\omega\otimes}{\longrightarrow}
\Om^{1/2}E^2_J\otimes E^2_J \to \Om^{1/2}E^2
$$
(recall that $E^2_J$ is our shorthand for the  weightless version of
$\Lambda^+$). In other words, $2$ of the $5$ components of $W^+$ are
automatically zero in the Hermitian case.  Accordingly the
linearization of the  ASD  equations is given by a second-order
operator
$$
D_{c,J}:\ci(\Om^{1/4}E^1_J) \to \ci(\Om^{3/4}E^2_J) 
$$
where $(\omega\otimes D_{c,J} u)_0 = D_c(\omega\otimes u)$. If
$\omega$ satisfies the twistor equation, then $D_{c,J}$ simplifies to
a conformally invariant operator which we shall denote by $S$
\begin{equation}\label{Sdef}
S: \ci(X,\Om^{1/4}E^1_J) \to \ci(X,\Om^{3/4}E^2_J)
\end{equation}
which {\em depends only upon $c$} and not upon $\omega$.  This
$S$-operator was studied in detail in \cite{LebS}; it has the form
$d^+d^* + \Phi$ where $\Phi:\Om^{1/4}E^1_J \to \Om^{3/4}E^2_J$ is
defined by regarding the trace-free part of the Ricci tensor as a
section of $\mbox{Hom}(\Om^{1/4}E^1_J,\Om^{3/4}E^2_J)$. The operator
$D_{c,J}$ or $S$ has the same role, in the study of Hermitian-ASD
conformal structures, as the deformation complex \eqref{deform-cx} for
ASD conformal structures.  In particular if we put
\begin{equation} \label{hjc}
H^1_{c,J}(X) = \ker(S),\; H^2_{c,J}(X) = \coker S = \ker S^*
\end{equation}
then $H^1_{c,J}$ is formally the tangent space to the moduli space of
hermitian-ASD deformations of $c$, while $H^2_{c,J}$ is an obstruction
space.

\subsubsection{Nonlinear terms}

The analysis of nonlinear terms in the hermitian-ASD equations is
the same as for the full ASD equations and leads to the following,
which is the analogue of Proposition~\ref{asdsum}:
\begin{propn} \label{hasdsum}
Let $(X,c,J)$ be a conformal hermitian surface, let $\omega$ be the
(conformal) fundamental $2$-form, 
and let
\begin{equation} \label{wtdjbunds}
E^1_J = \Om^{-1/2}\Lambda^-,\;\;E^2_J = \Om^{-1/2}\Lambda^+.
\end{equation}
Then
 there exists a $\ci$ map  $F_J: B_J \to C^\infty(X,\Om^{3/4}E_J^2)$, 
where $B_J$ is an open neighbourhood of the zero-section in
$\ci(X,\Om^{1/4}E^1)$, with an expansion of
the form
\begin{equation}\label{Fjdef}
  F_J(u) = W^+_0 + D_{c,J}u + \varepsilon_1(1+u,u\otimes\nabla\nabla u) +
    \varepsilon_2(1+u,\nabla u\otimes \nabla u) +
    \varepsilon_3(1+u,u\otimes u).
\end{equation}
such that $F_J(u) =0$ if and only if $c_{\omega\otimes u}$ is
Hermitian-ASD. In \eqref{Fjdef},
$W^+_0$ is the self-dual Weyl tensor of $c$, $D_{c,J}$ is the
linear operator mentioned above and the $\ve_j$ are nonlinear
terms which conform to Convention~\ref{conv}. Finally if $\omega$ satisfies
the twistor equation, so that $(c,\omega)$ is conformal to scalar-flat
K\"ahler, then $D_{c,J}$ can be replaced by the operator $S$ of \eqref{Sdef}.
\end{propn}

This completes our study of the PDE aspects of the two problems of
interest in this paper.  We have seen that when regarded as
`perturbation problems', both have elliptic linearizations (in the
case of the full ASD equations it is necessary to divide by the action
of the diffeomorphism group). We are going to apply
Propositions~\ref{asdsum} and \ref{hasdsum} on the manifold $X_\rho$ to expand the ASD
equations about the approximate solution $g_\rho$ constructed in
\S\ref{glueb}.  In the next section the behaviour of the
linearizations $D_c$ and $D_{c,J}$ will be studied, leading to an
application of the implicit function theorem to solve the equations in \S\ref{s5}.

\section{Linear theory on $X_\rho$}
\setcounter{equation}{0}
\label{linear}

This section is devoted to a discussion of the linear aspects of our
gluing problems. The first part includes the `main estimate' 
(Proposition~\ref{prmest}) which is one of the key technical results
needed for all our gluing theorems. This proposition states, roughly
speaking, that whenever  a fully elliptic operator $P_\rho$ on
$X_\rho$ is obtained by gluing fully elliptic operators $P_j$ on $X_j$
then $\ker P_\rho$ and $\coker P_\rho$ are well approximated by 
$\ker P_1 \oplus \ker P_2$ and $\coker P_1 \oplus \coker P_2$, 
and that $P_\rho$ induces a uniformly bounded isomorphism between 
suitable complementary subspaces.

In \S\ref{applic} we discuss the linearized operator that arises in the
deformation complex and note that it is never fully elliptic on a 
$b$-manifold. This entails a discussion of the linearization over 
cylinders which is taken further than is strictly needed for most of 
the gluing theorems. 

In \S\ref{compty} we consider $b$-manifolds that are conformal blow-ups and
blow-downs of compact or WALE spaces and compare 
the `Hodge' version of the deformation complex
for a $b$-manifold with the analogous definitions in terms of the compact 
or WALE models. Here conformal invariance is  an
essential tool.

Finally  in \S\ref{lthyhasd} we give comparison 
results for the operator $S$ 
that arises in the hermitian-ASD problem. 

\subsection{Stretching the neck}
\label{stretch} In this section we consider the
behaviour of linear elliptic operators  on $X_\rho$ as $\rho\to\infty$,
under the assumption that the corresponding operators over $X_0$,
$X_1$ and $X_2$ are {\em fully elliptic}. The main points of the
argument are based closely on Floer's paper, but with several
simplifications and generalizations.  Because the analysis presented
here is also needed in other geometric applications \cite{ko,KS}, 
we work here in some generality. 

\subsubsection{Notation}
\label{notn4}
The geometric set-up and notation used will be as in \S\ref{glueb}.  
In particular $X_1$ and $X_2$ will be $b$-manifolds and  $Y$ a piece
of $\p X_j$ at which the gluing is taking place.

Now suppose that for $j=1,2$,
$$
P_j: C^\infty(X_j,E_j) \to C^\infty(X_j, F_j)
$$
are fully elliptic phg $b$-differential operators of order
$m$. Suppose that over the cylindrical subset $U_j$ of \S\ref{data} we
have identifications $E_j = \pi^*E_0$, $F_j = \pi^*F_0$, where
$\pi:U_j \to Y$ is the obvious projection. Suppose further that 
with these identifications we have
$P_j = p_j(x_j,y;x_j\p_{x_j},\p_y)$ and 
the indicial operators agree in the sense 
$$
P_0:=p_1(0,y;\p_t,\p_y) = p_2(0,y;-\p_t,\p_y).
$$
(The sign in $p_2$ is to take care of the sign convention in the
definition of the $t_j$.) 
Then we can glue the $P_j$ across $\pm t_j =\rho$ to 
obtain bundles $E_\rho$, $F_\rho$ and a fully elliptic $b$-operator
$$
P_\rho :C^\infty(X_\rho, E_\rho) \to C^\infty(X_\rho, F_\rho)
$$
in exactly the same way that the metric $g_\rho$ was constructed from
the $g_j$ in \S\ref{gbr}. We have analogous estimates to those of
\S\ref{firstp}: 
\begin{equation} \label{pest}
\sup_{|t|\leq T}|P_j - P_0|= O(e^{-\eta\rho}),\; \sup_{|t|\leq T}|P_\rho - P_0|, 
 = O(e^{-\eta \rho})\mbox{ as }\rho \to \infty
\end{equation}
for any fixed $T$, where $\eta>0$ is some constant. Here we have
written $|P_j - P_0|$ for the sum of the moduli of the coefficients of
$P_j - P_0$.  Note that in our application, the $P_j$ arise as
operators canonically associated to some geometric data (a metric or
conformal structure). Then by cutting off and gluing these data and
{\em then} constructing the corresponding $P$-operator one will get a
slightly different operator from $P_\rho$, constructed by gluing the
$P_j$ directly. Since in
either case we shall have estimates like those of \eqref{pest}, this
difference is not important, and will be ignored in the sequel.

\subsubsection{Definition of asymptotic kernels and cokernels}

According to Proposition~\ref{tinbt}, 
$$
P_0 : L^p_m(X_0,E_0) \to L^p(X_0,F_0)
$$
is an isomorphism and by Theorem~\ref{frpk}, 
the $P_j$ are Fredholm in $L^p$ for every
$p$, with index independent of $p$. Let the $L^p$ null-space of $P_j$
be denoted by $N_j$ and the $L^p$ null-space of the $L^2$ adjoint
$P_j^*$ be denoted by $M_j$. By part (ii) of Theorem~\ref{frpk}, $N_j$
and $M_j$ consist of phg  sections. In particular, these sections are 
exponentially decreasing as $t_1 \to \infty$ or $t_2 \to
-\infty$, and the same is true of all derivatives of these sections.
Because of this,
$$
N_{1,\rho} = \beta(t_1 -\rho/2)N_1,\;\;M_{1,\rho} = \beta(t_1
-\rho/2)M_1
$$
and 
$$
N_{2,\rho} = \beta(-t_2 +\rho/2)N_2,\;\;M_{2,\rho} = \beta(-t_2
+\rho/2)M_2
$$
obtained by cutting off at distance $\rho/2$ will be very good
approximations if $\rho$ is large. Morever, we can regard
$N_{j,\rho}$ as a subspace of $U_\rho = L^p_m(X_\rho,E_\rho)$ and
$M_{j,\rho}$ as a subspace of $V_\rho = L^p(X_\rho,E_\rho)$. Define
$U'_\rho = N_{1,\rho}\oplus N_{2,\rho}$ to be the {\em asymptotic
kernel} of $P_\rho$ and $V'_\rho = M_{1,\rho}\oplus M_{2,\rho}$ to be
its {\em asymptotic cokernel}. Finally denote by $U''_\rho$
(resp.\ $V''_\rho$) the $L^2$-orthogonal complement of $U'_\rho$
(resp.\ $V'_\rho$) in $U_\rho$ (resp.\ $V_\rho$).

\subsubsection{The main estimate}

\begin{propn} \label{prmest}
There exists $\rho_* >0$ such that for all $\rho > 
\rho_*$, the induced map $P''_\rho : U''_\rho \to V''_\rho$ is an 
isomorphism and the operator norm of  $G_\rho = [P''_\rho]^{-1}$ is 
bounded independent of $\rho$.
\end{propn}

\noindent{\em Proof.} We shall prove that there exists $\rho_*$ and a 
constant $\varepsilon$ such that if $\rho > \rho_*$,
then
\begin{equation}\label{mest}
\|P''_\rho u\| \geq \varepsilon \|u\|\mbox{ for all }u\in U_\rho''.
\end{equation}
Such an estimate shows that $P_\rho$ is injective, with a uniform 
estimate for a left-inverse $V''_\rho \to U''_\rho$. Our set-up is 
symmetric under adjoints, however, so by the analogous estimate with
$P_\rho$ replaced by $P_\rho^*$, we see that
$P''_\rho$ is surjective, 
and the left-inverse is a true inverse, with norm bounded independent
of $\rho$.

The proof of \eqref{mest} goes by contradiction. If it fails, there
exists a sequence $\rho_n \to \infty$ and $u_n \in U''_n$ such 
that 
\begin{equation} \label{cont1}
\|u_n\| = 1
\end{equation}
but
\begin{equation} \label{cont2}
\|P''_n u_n\| \to 0
\end{equation}
as $n \to \infty$. (Here we begin to use obvious notational
simplifications, replacing $\rho_n$ by $n$ wherever this is unambiguous.)
From the definition of $P''_\rho$, there exist 
further sequences $v_n^{(j)}\in V'_{j,n}$ such that
\begin{equation} \label{cont3}
\|P_n u_n - v_n^{(1)} -v_n^{(2)}\| \to 0
\end{equation}
as $n\to \infty$. The main step in the proof is contained in the 
following:

\begin{lem} \label{lemming}
There exists a subsequence of $u_n$ (which by abuse of 
notation we continue to denote by $u_n$), such that
$$
\|u_n\|_{L^p_m(|t| \leq 2)} \to 0\mbox{ as }n\to \infty.
$$
\end{lem}

(We shall see from the proof that $2$ could be replaced by any larger
positive real number.)

\noindent{\em Proof of Lemma~\ref{lemming}} Multiply $u_n$ by a 
bump-function
so that it is cut off to zero at $t = \pm \rho_n/2-1$ Denote by
$u^{(0)}_n$ the resulting section over $X_0$. Clearly the $L^p_m$-norm
of $u^{(0)}_n$ is uniformly bounded as $n\to \infty$, so by passing to
a subsequence we may assume that $$ u^{(0)}_n \to u^{(0)}\mbox{ weakly
in }L^p_m.  $$ Over any fixed compact $K\subset X_0$, we have by \eqref{cont3}
$$
\|P_0u_n^{(0)}\|_{L^p(K)} \leq \|(P_0 - P_n)u_n\|_{L^p(K)} + 
\|P_n u_n\|_{L^p(K)} \to 0\mbox{ as }n\to \infty
$$ 
since the support of $v_n^{(1)} + v_n^{(2)}$ does not meet $K$.
Now the seminorm $u \mapsto \|P_0 u\|_{L^p}$
is continuous with respect to the $L^p_m$-norm, so by
lower-semicontinuity of weak limits, $\|P_0 u^{(0)}\|_{L^p(K)}=0$.
Since this is true for every compact set $K$, $P_0u^{(0)}=0$ over
$X_0$ and hence $u^{(0)}=0$ because $P^{(0)}$ is an isomorphism in $L^p$.

To complete the proof of the lemma, recall that weak convergence in $L^p_m(K)$ 
implies strong convergence in $L^p(K)$ if $K$ is compact.  The Lemma 
now follows by taking $K = \{|t| \leq 2\}$ and applying the elliptic 
estimate
$$
\|u_n\|_{L^p_m(K)} \leq C(\|P_nu_n\|_{L^p(K)}+ \|u_n\|_{L^p(K)}).
$$

{\em Proof of Proposition~\ref{prmest}}.  Replace $u_n$ by the subsequence 
in Lemma~\ref{lemming} and construct sequences $u_n^{(j)}$ over $X_j$ by 
cutting off $u_n$ in the intervals $\rho_n-2 < t_1 < \rho_n -1$ and
$1-\rho_n < t_2 < 2-\rho_n$. More precisely, we have 
$u_n^{(j)} = \beta_n^{(j)} u_n$, say, where $\beta_n^{(j)}$ is a 
suitable translation of a standard cut-off function.  We 
have
$$
P_j u^{(j)}_n = \beta_n^{(j)} P_n u_n + [P_j,\beta_n^{(j)}] u_n
$$
since $P_n = P_j$ where $\beta_n^{(j)}\not=0$. Therefore
$$
\|P_j u^{(j)}_n - v_n^{(j)}\| \leq \|P_n u_n - v_n^{(1)} - 
v_n^{(2)}\| + \|[P_j,\beta_n^{(j)}] u_n\|
$$
and the first term here tends to zero by \eqref{cont3}, the second by 
Lemma~\ref{lemming}, for $[P_j,\beta_n^{(j)}]$ is a differential
operator of order $\leq m$ which vanishes outside of $\{|t|\leq 2\}$.
 Since $P_j u^{(j)}_n$ and $v_n^{(j)}$ lie in
complementary subspaces, it follows that
$$
\|P_j u^{(j)}_n\| \to 0 \mbox{ and }\|v^{(j)}_n\| \to 0.
$$
On the other hand $u^{(j)}_n$ lies in a subspace on which $P_j$ is 
injective, so $u^{(j)}_n \to 0$ in $L^p_m$ as $n\to \infty$.

Finally we combine these estimates with \eqref{cont1} to obtain
$$
1 = \|u_n\|  \leq \|\beta^{(1)}_n u_n \| + \|\beta^{(2)}_n u_n \|
+ \|(1-\beta^{(1)}_n -\beta^{(2)}_n)u_n\| \to 0.
$$
This contradiction completes the proof of Proposition~\ref{prmest}.

\subsection{Application to ASD problem}
\label{applic}
Let $X$ be a compact oriented $4$-manifold with a $b$-metric
$g$. Let
\begin{equation}
\label{frmadj} \DD=(D_g,L^*_g):\ci(X, E^1) \to \ci(X,E^2)\oplus \ci(X, E^0)
\end{equation}
so that if $\p X=\emptyset$ and $g$ were conformally ASD we should have
$$
H^1_c = \ker \DD,  H^2_{c} \oplus H^0_{c} = \coker \DD
$$
the direct-sum decomposition of the cokernel corresponding to the
direct-sum decomposition in \eqref{frmadj}.  We wish to apply the main
estimate to $\DD_\rho$ obtained by gluing $\DD_j$. Therefore we need
to check whether $\DD_g$ is fully elliptic when $g$ is a
$b$-metric. For this we must examine the indicial operator $I(\DD)$ or
equivalently the operator $\DD_0$ corresponding to the cylinder $X_0 =
Y\times \RE$, where $g_0=h(y) +
dt^2$ is a product metric.  We assume that 
 $g_0$ is conformally ASD, hence
conformally flat (because $t$-invariant). Since $dt^2 + h(y)$ is
conformally flat if and only if $h$ is a metric of constant
curvature we assume this from now on. Then the
section $(dt\otimes dt)_0$ of $E^1$ is parallel so $\DD_0(dt\otimes
dt)_0 = (\Phi_0\cdot (dt\otimes dt)_0$ where $\Phi_0$ is the 
trace-free part of the Ricci tensor of $g_0$ (which is
itself a multiple of $(dt\otimes dt)_0$) and $(\cdots)_0$ denotes the
projection into $S^2_0\Lambda^+$. Hence this term is $0$ and 
it follows that $\DD_g$ is {\em never} fully elliptic.

We can however apply the conjugation trick of \S\ref{cjw} to obtain
nearby operators that {\em are} fully elliptic. To this end, 
denote by $\DD^{\pm}_g$ the operator
\begin{equation} \label{dpmdef}
\DD_g: x^{\pm\delta}L^2_2(X,E^1) \to x^{\pm\delta}(L^2(X,E^2)\oplus
L_1^2(X,E^0))
\end{equation}
and set
$$
\ker \DD_g^{\pm} = H^{1,\pm}_c,\;\;\coker \DD_g^{\pm} = H^{2,\pm}_c\oplus
H^{0,\pm}_c, 
$$
the direct-sum decomposition of the cokernel corresponding, as before,
to that of the target of $\DD$.  Here we assume
\begin{equation} \label{drange}
\delta >0, \{0<|\Im\lambda|\leq \delta\}\cap\spb(\DD_g) = \emptyset.
\end{equation}
It follows that $\DD^\pm$ are Fredholm and
that the $H^{j,\pm}$ are independent of the choice of
Sobolev spaces in \eqref{dpmdef} and of the choice of $\delta$ satisfying \eqref{drange}.

From the definitions, 
$$
H^{1,+}_g\subset H^{1,-}_g,\; H^{2,-}_g\subset H^{2,+}_g,\;
H^{0,-}_g\subset H^{0,+}_g.
$$
We are now in a situation in which Proposition~\ref{prmest} can be
applied.  We take this up in \S5 and the reader who wants to see at 
once how this leads to gluing theorems for ASD conformal structures 
can now start reading there.

\vspace{5pt}
The remainder of this section is however taken up with a more
detailed discussion of the indicial operator $\DD_0$.

\subsubsection{$\DD$ on a cylinder} \label{dcyl}
With  $(X_0,g_0)$ the
cylinder $Y\times \RE$ with product metric $g_0=h(y) + dt^2$, 
where $h$ is
a metric of constant curvature, we wish to study the non-vanishing
$H_c^{j,\pm}(Y\times \RE)$, viz.\
$H_c^{0,+}(Y\times \RE)$, $H_c^{1,-}(Y\times \RE)$, 
$H_c^{2,+}(Y\times \RE)$.  
By definition, these are respectively the parts of the the
null-spaces of $L$, $\DD$, $D^*$, that consist of $t$-invariant or
$t$-periodic sections. An index theorem gives 
$$
\dim H_c^{0,+}(Y\times \RE)-\dim H_c^{1,-}(Y\times \RE) +
\dim H_c^{2,+}(Y\times \RE) = 0
$$
for these spaces can be identified with cohomology groups for an
associated elliptic complex over $Y$.  Furthermore, it is not hard to
see that $H_c^{0,+}(Y,\times \RE) = \RE \oplus \mbox{Isom}(Y)$, where
the summand $\RE$ corresponds to translations of the cylinder.

\subsubsection{Geometric interpretation of $H_c^{1,-}(Y\times \RE)$} 
\label{gi} If $H_c^{1,-}(Y\times \RE)\not=0$, then
there exists a $t$-periodic or $t$-invariant solution of $\DD_0 u =
0$. In particular, we can consider $u$ to be a solution on $Y\times
S^1$, where the length of the circle is determined by the period of
the solution. Thus $u$ gives an infinitesimal conformally flat deformation of the
product conformally flat structure on $Y\times S^1$. Now there is an
obvious space of deformations of this structure: constant-curvature deformations of the
constant-curvature metric on $Y$ plus isometries of $Y$ used to change
the way in which the two copies of $Y$ are identified when making
$Y\times S^1$ from $Y\times I$ ($I$ is an interval). Such deformations
correspond to $t$-invariant solutions $u$, and in the cases which are
understood, when the curvature of $Y$ is non-negative, all solutions are
of this form.

\subsubsection{Positive curvature: $Y= S^3/\Ga$} According to the
previous paragraph, the $t$-invariant part of $H^{1,-}$ contains 
$\RE \oplus  {\mathfrak so}_4^\Ga$. The $\RE$ summand corresponds to
homotheties of $h(y)$ (or to $(dt\otimes dt)_0$) while
${\mathfrak so}_4^\Ga$ is the space of infinitesimal 
isometries of $S^3/\Ga$. Recall that homotheties are the only
constant-curvature deformations of the round metric on $S^3$.
Floer has shown, moreover, that this is the whole of
$H^{1,-}_0(S^3\times \RE)$ \cite[\S 5]{floer-conformal}.
Summing up,
$$
H^{0,+}_c(S^3/\Ga\times \RE) = H^{1,-}_c(S^3/\Ga\times\RE) =
\RE \oplus  {\mathfrak so}_4^\Ga,\;\;
H^{2,+}_0(S^3/\Ga\times \RE) = 0
$$
by the index formula in \S\ref{dcyl}.\label{pscrv}

\subsubsection{The case $Y = T^3$} 
The metric on $Y\times \RE$ is in
this case flat and it follows from
$\DD u =0$  that $\Delta^2u=0$. If $u$ is also bounded, then
$u$ must be
parallel. In particular, every solution is $t$-invariant and
$H^{1,-}_c(T^3\times \RE) = \RE^9$. (This also agrees with the geometric
interpretation of \S\ref{gi}.)  Similarly we have
$$
H^{0,+}_c(T^3\times \RE) = \RE^4,\;
H^{2,+}_c(T^3\times \RE) = \RE^5.
$$

\subsubsection{The case $Y$ is hyperbolic} In this case,
\S\ref{gi} would predict that $H^{1,-}_0(Y\times \RE) = \RE$ because
of the rigidity of hyperbolic structures. In fact, by following
Floer's calculations through one can verify that if $Y$ is a
hyperbolic rational homology $3$-sphere, then the $t$-invariant part
of $H^{1,-}_0(Y\times \RE)$ is indeed $1$-dimensional.  
We have not been able to eliminate
the possibility of periodic solutions in this case, but 
make the following
\begin{conjecture} If $Y$ is a hyperbolic rational homology
$3$-sphere, then $H^{1,-}_0(Y\times \RE) = \RE$.
\end{conjecture}
In any case, we have $H^{0,+}(Y\times \RE) = \RE$. If the conjecture
is true, then we should also have $H^{2,+}_0(Y\times \RE)=0$.

This conludes our
digression on the properties of $\DD_0$.

\subsection{Comparison theory}
\label{compty}
We now turn to $b$-manifolds that arise by conformal blow-up or
blow-down (cf.\  \S\ref{confb}, \S\ref{walex}). The main results
of this section are the following `comparison theorems':

\begin{thm} Let $(X,g)$ be the conformal blow-up at the point $0$ of
the compact conformally ASD orbifold $(\bX,\bg)$. Then the framed cohomology
groups $H^*_{\bc}(\bX,0)$ of the  framed deformation complex
$$
C^\infty_{(0)}(\bX,\Om^{-1/4}E^0) \to
C^\infty(\bX,E^1) \to C^\infty(\bX,\Om^{1/2}E^2),
$$
where 
$$
\ci_{(0)}(\bX,V) = \{v\in \ci(\bX,V): v(0) =0\},
$$
agree with the kernel and cokernel of $\DD^+_g$:
$$
H^j_{\bc}(\bX,0) = H^{j,+}_g(X), \mbox{ for }j=0,1,2.
$$
Moreover, 
$$
H^{2,+}_c(X) = H^2_{\bc}(\bX) =
\ker(\bD^*:\ci(\bX,\Om^{1/2}E^2) \to
C^\infty(\bX,\Om E^1) ).
$$
\label{comp1}
\end{thm}

This result was stated
and partly proved by Floer \cite[\S 3]{floer-conformal}. The following
is the analogue for comparison with a WALE space:

\begin{thm} Let $(X,g)$ be the conformal blow-down at $\infty$ of
the  conformally ASD, WALE space $(\hX,\hg)$. Assume that $\hg$ is
such that $g$ is a smooth $b$ metric. Then the cohomology
groups $H^*_{\hc}(\hX)$ of the  complex
$$
C^\infty(\hX,\Om^{-1/4}E^0)\cap O(1) \to
C^\infty(\hX,E^1)\cap O(R^{-1}) \to
C^\infty(\hX,\Om^{1/2}E^2)\cap O(R^{-3}),
$$
where 
$$
\ci(\hX,V)\cap O(R^{-a}) = \{v\in \ci(\hX,V): |\nabla^j v | =
O(R^{-a-j}),\mbox{as $R\to\infty$, for all }j\}
$$
agree with the kernel and cokernel of $\DD^-_g$:
$$
H^j_{\hc}(\hX) = H^{j,-}_g(X), \mbox{ for }j=0,1,2.
$$
Moreover,
$$
H^{2,-}_c(X) = H^2_{\hc}(\hX) =
\ker(\hD^*:\ci(\hX,\Om^{1/2}E^2)\cap O(R^{-4}) \to
C^\infty(\hX,\Om E^1) ).
$$
\label{comp2}
\end{thm}

Note that the decay conditions used to define the $H^j_{\hc}$ are the
natural ones on a WALE space, In particular,
a differential operator of order
$m$ canonically associated with $\hg$ will automatically map 
$O(R^{-a})$-sections to $O(R^{-a-m})$-sections.

\begin{proof} Recall first the relation between the cylindrical
metric, and the euclidean metrics $g_0$ and $g_\infty$ (near $0$ and
$\infty$, respectively):
\begin{equation} \label{ceq1}
dt^2 + d\omega^2 = r^{-2}(dr^2 + r^2 d\omega^2) = R^{-2}(dR^2 + R^2 d\omega^2) 
\end{equation}
where
\begin{equation} \label{ceq2}
t = e^{-r} = e^{-x},\; t = e^R,\;\; R = r^{-1}.
\end{equation}
Here we are assuming that $0$ is a smooth point, or if not, we pass to
a uniformizing chart centred at $0$. It will be clear that
$\Gamma$-equivariance is preserved throughout, so we can afford to 
ignore singularities from now on.
Denoting the pointwise norms that correspond to the different choices
of conformal gauge in \eqref{ceq1} by $|\cdot|_{b}$,
$|\cdot|_0$ and $|\cdot |_\infty$ respectively, we have by
\eqref{metch} and \eqref{mtdef},
\begin{quote}

(i) If $\xi \in \Om^{-1/4}E^0$, then 
$|\xi|_b = r^{-1}|\xi|_0 = R^{-1}|\xi|_\infty$;

(ii) If $h \in E^1$, then 
$|h|_b = |h|_0 = |h|_\infty$;

(iii) If $\psi \in \Om^{1/2}E^2$, then 
$|\psi|_b = r^{2}|\psi|_0 = R^{2}|\psi|_\infty$.
\end{quote}
We compare first $H^0$ and $H^2$. By the weighted Fredholm
alternative, we have
$$
H^{0,\pm}_c(X)  = \ker(x^{\mp\delta} L^p_k(X,\Om^{1/4}E^0) 
\stackrel{L}{\longrightarrow}
x^{\mp\delta} L^p_{k-1}(X,E^1))
$$ and
$$
H^{2,\pm}_c(X)  = \ker(x^{\mp\delta} L^p_k(X,\Om^{1/2}E^2) 
\stackrel{D^*}{\longrightarrow}
x^{\mp\delta} L^p_{k-1}(X,\Om E^1)),
$$
the formal adjoint being taken with respect to the $b$-metric $g$.
Therefore by conformal invariance of these differential operators and
the rescaling formulae (i) and (iii) above,
$$
H^{0,\pm} = \{\bxi\in \ci(\bX\backslash \{0\},\Om^{-1/4}E^0): \bL\bxi =
0,\; |\xi|_0 = O(r^{1\mp\delta})\mbox{ as }r\to 0\},
$$
and
$$
H^{2,\pm} = \{\bpsi\in \ci(\bX\backslash \{0\},\Om^{1/2}E^2): \bD^*\bpsi =
0,\; |\psi|_0 = O(r^{-2\mp\delta})\mbox{ as }r\to 0\}.
$$
We argue next that the singularity at $0$ is removable in all cases
except $H^{2,+}$.  First note that in all cases, $\bxi$ and $\bpsi$
admit extensions as distributional sections to the whole of
$\bX$ \cite[\S3.2]{h1}. Denoting these extensions by the
same symbols, we obtain equations of the form $\bL\bxi = \alpha$,
$\bD^*\bpsi = \beta$, where $\alpha$ and $\beta$ are linear
combinations of derivatives of the Dirac distribution $\delta_0$.
Thus $\alpha$ and $\beta$ are sums of terms homogeneous of degree
$-4-m$ for $m=0,1,\ldots$ and this is not compatible with the orders
of $\bL$ and $\bD^*$ and the growth conditions on $\bxi$ and $\bpsi$
that appear in $H^{0,\pm}$ and $H^{2,-}$. So in these three cases,
$\alpha=0$ and  $\beta=0$, and since $\bD^*$ and $\bL$ are
overdetermined elliptic the solutions are actually smooth. Hence we
obtain
\begin{eqnarray}
H^{0,+} & = &\{\bxi\in \ci(\bX,\Om^{-1/4}E^0): \bL\bxi =
0,\; \bxi(0) = 0\} \\
H^{0,-} & = &\{\bxi\in \ci(\bX,\Om^{-1/4}E^0): \bL\bxi =
0,\; \bxi(0) = 0, \bdel\bxi(0)=0\} \\
 &= & \{\hxi\in \ci(\hX,\Om^{-1/4}E^0): \hL\hxi =
0\} \cap O(1),
\end{eqnarray}
the last following from the formula (i) relating the lengths of an
element of $\Om^{1/4}E^0$ in the compact and WALE models.  We also
have
\begin{eqnarray}
H^{2,-} & = &\{\bpsi\in \ci(\bX,\Om^{1/2}E^2): \bD^*\bpsi = 0 \} \\
&=&  \{\hpsi\in \ci(\hX,\Om^{1/2}E^2): \hD^*\hpsi = 0\} \cap  O(R^{-4}).
\end{eqnarray}
In order to overcome the problem encountered with $H^{2,+}$ we note a
result of Biquard \cite{biquard} which provides 
exact comparisons of $b$-Sobolev spaces of $X$ with
ordinary Sobolev spaces on $\bX$ for a good choice of weight $\delta$
and exponent $p$. The version we need states that the obvious map on
compactly supported functions extends to an isomorphism
$$
\lambda: x^\delta L^p_2(X,E^1) \simeq L^p_2(X,0;E^1)
:=\{h\in L^p_2(\bX,E^1): h(0) =0\}
$$
provided that 
\begin{equation} \label{dpcond}
0 < \delta = 2 -4/p < 1.
\end{equation}
  Note that the vanishing
condition makes sense because $L^p_2 \subset C^0$ for $2 -4/p>0$ in
$4$ dimensions and that Biquard's result for functions applies here
because $E^1$ is weightless (and so behaves conformally like the
trivial bundle). By direct calculation (in which the factor
$\Om^{1/2}$ is crucial) we see also that there is an isomorphism
$$
\mu: x^\delta L^p(X,\Om^{1/2}E^2) \simeq L^p(\bX;\Om^{1/2}E^2)
$$
if $\delta$ and $p$ are related as in \eqref{dpcond}.
Using conformal invariance of $D$ we have $D^+ = \mu^{-1}\bD\lambda$ and in
particular
$$
H^{2,+} = \coker(\bD:L^p_2(\bX,0;E^1) \to L^p(\bX,\Om^{1/2}E^2)) 
$$
It follows that  $H^2_{\bc}(\bX) = H^{2,+}_g$ if 
$$
\bD(L^p_2(\bX,E^1)) = \bD\{\bh \in L^p_2(\bX,E^1): \bh(0) = 0\}.
$$
To prove this, it is enough to show that given $h\in L^p_2(\bX,E^1)$,
there exists $\xi \in L^p_3(\bX, E^0)$ such that $h(0) - \bL\xi(0)
=0$. Working in normal coordinates $x^a$ near $0$, we have the formula
$$
(\bL\xi)_{ab} = \p_a\xi_b + \p_b\xi_a - [\p_c\xi^c]g_{ab}/2\mbox{ at }0
$$
just as in $\RE^4$ with the Euclidean metric. Then
$\bL(\beta\,h_{ab}x^a) (0) = 2h_{ab}$ if $h_{ab}$ is constant and
$\beta = 1$ in a neighbourhood of $0$.  This completes the
proofs of all statements pertaining to $H^0$ and $H^2$ in
Theorems~\ref{comp1} and \ref{comp2}.

\subsubsection{Comparison of $H^1$}  Fix $p$ and $\delta$ as in
\eqref{dpcond}. Then by Biquard's result and conformal invariance,
there is a natural map
$H^{1,+} \to H^1_{\bc}(\bX,0)$ given by mapping $h\in H^{1,+}_c$ to the
cohomology class  $[\bh]$ of the image of $h$ in
$L^p_2(\bX,0;E^1)$. We claim first that this map is injective. Indeed,
if $\bh = \bL\bxi$ for some $\bxi$ with $\bxi(0)=0$, then transferring
back to $X$ we obtain $\xi$, such
that $|\xi|_b = O(1)$ as  $t\to \infty$ and satisfying $L\xi = h$. But
$h\in H^{1,+}_c$ implies that $L^*L\xi =0$ and since 
$h$ is in $L^2$, we may integrate by parts, getting $h=L\xi=0$. This
establishes the injectivity.

To prove surjectivity, we need to show that if $\bh$ represents a
class in $H^1_{\bc}$, then we can find $\bxi \in C^\infty(E^0)$, such
that $\bxi(0)=0$ and 
$$
L^*(L\xi +h) = 0
$$
By a previous argument we may
assume $\bh(0)=0$ so that $L^*h \in x^{\delta}L^p_2(X,E^1)$.  It is
straightforward to show
$$
L^*L(x^{\delta }L^p_k(X,E^0)) = L^*(x^{\delta }L^p_{k+1}))
$$
by adapting standard Hodge-theory arguments. In particular, we can
solve the equation and transfer to $\bX$, getting a section $\bxi$
which vanishes at $0$. This completes the comparison of $H^1$ for
conformal blow-ups, and so the proof of Theorem~\ref{comp1}. The
comparison of $H^1$ in Theorem~\ref{comp2} follows very similar lines,
using an analogue of Biquard's theorem to compare weighted Sobolev
spaces on $X$ with Sobolev spaces on $\hX$. The details are omitted.
\end{proof}

\subsection{Linear theory for the Hermitian-ASD problem}
\label{lthyhasd}
In \S\ref{hasd}, we saw that the operator 
$$
S: \ci(X,\Om^{1/4}E^1_J) \to \ci(X,\Om^{3/4}E^2_J),
$$
where $E^1_J$ and $E^2_J$ are the weightless versions of $\Lambda^-$
and $\Lambda^+$ respectively, controls the deformation theory of
conformally scalar-flat K\"ahler metrics.  In contrast to $\DD$, it
turns out that $S$ is fully elliptic on a $b$-manifold $X$ if each
component of $\p X$ is a spherical space-form.  This makes for
substantial simplifications, particularly for the comparison theorems
for this problem.

\begin{propn} Let $X_0 = S^3\times \RE$, $g_0 = h(y) + dt^2$,
where $h$ is the round metric on $S^3$, $S_0$ the $S$-operator
associated to $g_0$. Then if for some $u(y)$, we have
$S_0(u(y)e^{i\lambda t})=0$, 
it follows that $i\lambda \in \IZ \backslash 0$.
\label{sspec}
\end{propn}

 \begin{proof} Use the conformal isometry 
$S^3\times \RE \to \RE^4 \backslash 0$ given by $r= e^{t}$ and the conformal
invariance of $S$. Because of the conformal weights, we have if $u$
is a section of $\Om^{1/4}E^1_J$
$$
|u|_{\rm cylinder} = r|u|_{\rm Euclidean}
$$
so that a solution $\bS\bu =0$  in $\RE^4\backslash 0$, homogeneous of degree
$\lambda$,  translates into an exponential solution with factor
$e^{(\lambda +1)t}$ on the cylinder. In particular the constant
solution in $\RE^4$  gives rise to a solution that goes like $e^t$
along the cylinder.

We use a `removable singularities' argument like the one in the proof
of Theorem~\ref{comp1}. If $Su =0$ in $\RE^4\backslash 0$ and
$u$ has homogeneity $\lambda$ in $r$, then 
$Su$ is a distribution supported
at $\{0\}$ and homogeneous of degree $\lambda -2$.  If this
distribution vanishes at $0$, then by elliptic regularity, $u$ is
smooth near $0$ and hence $\lambda$ is a non-negative integer.  If the
distribution is non-vanishing, then
$Su$ must be a multiple of some derivative of  $\delta_0$, hence
$\lambda -2 = -4 - m$ for some integer $m\geq 0$. Hence $\lambda = -2
- m$ and so $\lambda = -1$ cannot occur.
\end{proof}

\subsubsection{Remark} On $T^3\times \RE$, however, the existence of
non-trivial parallel sections obstructs the full ellipticity of $S$ on
a $b$-manifold some of whose boundary components are tori.

\subsubsection{Comparison theory for $S$}  Consider now the situation
of \S\ref{compty}. We have:
\begin{thm} Let $(X,g)$ be a conformally ASD manifold with smooth $b$
metric, and suppose that $(X,g)$ is the conformal blow-up at $0$ in
$(\bX,\bg)$ and the conformal blow-down of $\infty\in (\hX,\hg)$.
Let $H^1_{c,J}(X)$, $H^2_{c,J}(X)$ be respectively the kernel and
cokernel of 
$$
S_g: L^p_k(X,\Om^{1/4}E^1_J) \to L^p_{k-2}(X,\Om^{1/4}E^2_J),
$$
let  $H^1_{\bc,J}(X)$, $H^2_{\bc,J}(X)$ be respectively the kernel and
cokernel of $\bS= S_{\bg}$, and let
$$
H^1_{\hc,J}(\hX) = \ker(\hS)\cap O(R^{-2}),
H^2_{\hc,J}(\hX) = \ker(\hS^*)\cap O(R^{-2}),
$$
(where the adjoint is taken relative to the WALE metric $\hg$). Then
$$
H^1_{c,J} = H^1_{\bc,J} = H^1_{\hc,J}\mbox{ and }
H^2_{c,J} = H^2_{\bc,J} = H^2_{\hc,J}.
$$
\label{comp3} 
\end{thm}
\begin{proof} If  $Su=0$, $u \in
L^p_k(X,E^1_J)$, then from Proposition~\ref{sspec}, 
$u$ has a phg expansion where the index set is
just the positive integers. Translating to $\bX$ and remembering the
conformal weight, we get $\bu$ which satisfies $\bS\bu =0$ and 
$|\bu|_{\bg} = O(1)$ near
$0$. Hence by elliptic regularity, the $L^p$-null space of $S$  on $X$ agrees
with the standard null-space of $\bS$ on $\bX$.  The argument for $H^2$ is the
same, for $H^2_{c,J}$ can be identified with the $L^p$ null-space of
$S^*$ and this is
conformally invariant as an operator
between bundles with the same
conformal weights as for $S$ (cf.\ \S\ref{sconfinv}).   

The argument
for comparison with $\hX$ is also closely analogous.\end{proof}

\setcounter{equation}{0}
\section{Nonlinear theory}
\label{s5}
We come now to the problem of finding an exactly conformally ASD 
$b$-metric on $X_\rho$ as a perturbation of the metric $g_\rho$ 
constructed in \S\ref{gbr}. More precisely, assume that $(X_j,g_j)$ in 
\S\ref{glueb} are conformally ASD $b$-manifolds so that $g_\rho$ is 
conformally ASD except in the damage zone $\{-1/2 \leq |t-1| \leq 
1/2\}$ near the middle of the neck. Using the conformal class $c_\rho$ 
of $g_\rho$ as the reference point $c$ in Proposition~\ref{asdsum}
we need to find a small ($\sup_{X_\rho}|h|<1$), smooth $h$ satisfying
\eqref{Fdef},
\begin{equation}\label{tbsv}
0 = F_\rho(h) 
= W^+_\rho + D_\rho h 
+ \varepsilon_1(1+h,h\otimes\nabla\nabla h)
+ \varepsilon_2(1+h,\nabla h\otimes\nabla h)
+ \varepsilon_3(1+h,h\otimes h)
\end{equation}
where $W^+_\rho := W^+(c_\rho)$, $D_\rho := D_{c_\rho}$.

The strategy is to solve this equation in a suitable Banach space by 
use of a version of the implicit function theorem (IFT). Then the 
solution will be proved to be $\ci$ and polyhomogeneous near $\p 
X_\rho$. For both parts of the argument it is useful to supplement
\eqref{tbsv} with the gauge-fixing condition
\begin{equation}\label{fixg}
    L^*_\rho h = 0
\end{equation}
(the $*$ denoting $L^2$ adjoint with respect to $g_\rho$). Then the 
linearization of the map 
\begin{equation} \label{fullmap}
    h \mapsto (F_\rho(h), L^*_\rho h)
\end{equation}
is equal to the elliptic operator $\DD_\rho$.

\subsection{Weak solution}

In the next few paragraphs we explain how to arrange matters so that 
this strategy can be successfully pursued. We shall need to choose 
Banach spaces so that $F_\rho$ extends to a $\ci$ map with uniform 
behaviour as $\rho \to \infty$.  This involves estimating certain 
Sobolev norms of the nonlinear terms in $F_\rho$. But we begin by 
adjusting $F_\rho$ so that the neck-stretching analysis of \S\ref{stretch} 
can be applied to its linearization.

\subsubsection{Introduction of weights} The main estimate,
Proposition~\ref{prmest} 
does not apply directly to $\DD_\rho$ because, as we saw at the start
of \S\ref{applic},  this operator does  not arise by 
gluing a pair of fully elliptic operators. Therefore we replace 
\eqref{fullmap} by 
\begin{equation}\label{wmap}
    u \mapsto (F^w_\rho(u), w^{-1}L^*_\rho(wu))\mbox{ where }
     F^w_\rho(u) :=w^{-1}F_\rho(wu), 
\end{equation}
where $w:=w_\rho$ is a suitable weight-function on $X_\rho$. 

\subsubsection{Definition of $w_\rho$} For $j=1,2$, let $w_j$ be equal
to a generically chosen positive power of a boundary defining function
for $\p X_j \backslash Y$. We assume $0\leq w_j \leq 1$ on $X_j$, with
$x_j=0$ only at $\p X_j \backslash Y$ and $w_j=1$ near $Y$.  Extend
the neck parameter $t$ smoothly to $X_\rho$ (and denote the extension
also by $t$) so that the range of $t$ is $[-\rho-1,\rho+1]$ and $t=
-\rho-1$ near $\p X_1 \backslash Y$ and $t= \rho+1$ near
 $\p X_2\backslash Y$. 
Now for $\delta$ satisfying the conditions of \eqref{drange}, we put
$$w_\rho = w_1w_2e^{-\delta(t + \rho + 1)}.  
$$ 
 If $\rho$ is fixed, we
have $0\leq w_\rho \leq 1$, with $w_\rho=0$ only at $\p X_\rho$;
$w_\rho$ is a power of a boundary defining function near $\p X_\rho$
and decreases exponentially along the neck.

The linearization of \eqref{wmap} is 
$P_\rho:= w_\rho^{-1}\DD_\rho w_\rho$ 
which is obtained by gluing the fully elliptic operators 
$$
P_1 = (w_1x_1^{\delta})^{-1}\DD_1(w_1x_1^{\delta}) $$ 
and 
$$ 
P_2 = (w_2x_2^{-\delta})^{-1}\DD_1(w_2x_2^{-\delta})
$$ 
as in \S\ref{notn4}.
According to the main estimate we can now invert $P_\rho$ in a
controlled way in Sobolev spaces over $X_\rho$.  The next task is to
choose Sobolev spaces such that \eqref{wmap} extends to a smooth map
between them.\label{wtdef}

\subsubsection{Choice of Sobolev space} Because $D_\rho$ is a
second-order operator, it is natural to take $u\in L^p_2(X_\rho,E^1)$
for some $p$. If $p>2$ we have the estimate
\begin{equation} \label{sob1}
\sup_{X_\rho} |u| \leq C\|u\|_{L^p_2(X_\rho)}
\end{equation}
and we can assume $C$ is independent of $\rho$. This uniformity of $C$
(and also \eqref{sob2}) follows from standard Sobolev embedding theorems for the complete
Riemannian manifolds $X_0$, $X_1$ and $X_2$, by a partition of unity
argument \cite[\S2.23]{aubin}.  We chose $w<1$ so that \eqref{sob1}
implies that
if the $L^p_2$-norm of $u$ is sufficiently small, then the
pointwise norm of $h=wu$ is everywhere $<1$ and $c_\rho(1+h)$ defines
a genuine ($C^0$) conformal structure. 

So $L^p_2$, for any $p>2$, will take care of the linear term in
\eqref{tbsv}. The nonlinearities in $F^w$ are much easier to control, however, if we
take $p>4$ so that \eqref{sob1} can be strengthened to
\begin{equation} \label{sob2}
\sup_{X_\rho}( |u|+|\nabla u|) \leq C\|u\|_{L^p_2(X_\rho)}
\end{equation}
for some other constant $C$ that is independent of $\rho$.
Therefore we fix $p>4$ and show next that $u\mapsto F^w(u)$ extends to
a smooth map from a neighbourhood of $0$ in $L^p_2(X_\rho,E^1)$ to
$L^p(X_\rho, E^2)$.

\subsubsection{Estimation of the nonlinearities}

From \eqref{tbsv} and the properties of the $\varepsilon_j$
(\S\ref{conv}) we obtain
\begin{equation}\label{fwdef}
(F_\rho^w(u),w^{-1}L^*(wu)) =
w^{-1}W_\rho^+  + P_\rho u +  Q(u)
\end{equation}
where 
$$
Q(u):=w\varepsilon_1(1+wu, u\otimes\nabla^w u\nabla^w u)
+ w\varepsilon_2(1+ wu, \nabla^w u\otimes \nabla^w u) 
+ w\varepsilon_3(1+wu, u\otimes u)
$$
and $\nabla^w u= w^{-1}\nabla (wu)$. (In the interests of legibility
we have dropped the notational
dependence of the weight $w$ upon $\rho$. There seems little danger of
confusion.) We can now state the main result:

\begin{propn} For any fixed $p>4$, the map 
$$
f_\rho:u \mapsto (F_\rho^w(u),w^{-1}L^*_\rho(wu))
$$
extends to a smooth map from the ball $B_\rho$ of radius $r$ in
$U_\rho = L^p_2(X_\rho,E^1)$ to $V_\rho$
where
$$V_\rho = L^p(X_\rho, E^2)\oplus L_1^p(X_\rho, E^0).$$
Moreover for a suitable choice of the parameters in the definition of
$w$, 
$f_\rho$ has the
following properties:\begin{itemize} 
\item[(i)] $\|f_\rho\| \to 0$ as $\rho \to \infty$;

\item[(ii)] There exist decompositions $U_\rho = U'_\rho\oplus
U''_\rho$,  $V_\rho = V'_\rho\oplus V''_\rho$, where $U'_\rho$ and
$V'_\rho$ are finite-dimensional, $U''_\rho$ and $V''_\rho$ are closed
and the map $P''_\rho:U''_\rho \to V''_\rho$ induced by $P_\rho$ 
has a uniformly bounded inverse $G_\rho:V''_\rho \to U''_\rho$.

\item[(iii)] The nonlinearity $Q(u)$
 satisfies
\begin{equation}
\|Q(u) - Q(v)\| \leq C(\|u\| + \|v\|) \|u-v\|
\end{equation}
for every $u$, $v$ in $B_\rho$, where $C$ is bounded independent of $\rho$.
\end{itemize}
\label{nlprop}
\end{propn}

\begin{proof} (i) Note first by the discussion in \S\ref{firstp} that
the $L^p$-norm of $W^+_\rho$ is $O(e^{-\eta\rho})$ for some $\eta
>0$. Taking $0< \delta < \eta$, 
$\|f_\rho(0)\| = O(e^{-(\eta-\delta)\rho}) \to 0$ as $\rho \to
\infty$. To check that the map is smooth, it evidently suffices to
show that the nonlinear terms define a smooth map. Note first by the
properties of $w$ that $\|\nabla^w u\|_{L^p_k} \leq
C\|u\|_{L^p_{k+1}}$ where $C$ is independent of $\rho$. Hence from
\eqref{sob2} we have
$$\|u\nabla^w\nabla^w u\|_{L^p} \leq (\sup|u|)\|\nabla^w\nabla^w
u\|_{L^p}
\leq C\|u\|^2_{L^p_2}.
$$
Similarly
$$
\|\nabla^w u\otimes\nabla^w u\|_{L^p} \leq C\sup (|u|+|\nabla
u|)\|u\|_{L^p_1}
\leq C\sup (|u|+|\nabla u|)\|u\|_{L^p_1}
\leq C\|u\|^2_{L^p_2}
$$
and
$$
\|u\otimes u\|_{L^p} \leq C\|u\|^2_{L^p_2}.
$$
(Here $C$ is a generic constant bounded independent of $\rho$ but
possibly varying from line to line.) This is not quite enough because
the $\varepsilon_j$ also have a real-analytic dependence on the
$0$-jet of $u$.  However, this is convergent for all $u$ with $\sup|u|
< 1$ and by multiplication properties of elements of $L^p_2$ with
$p>4$, these extend to define smooth maps from a fixed ball $B\subset
U_\rho$ into $V_\rho$. Combining these observations with the previous
estimates for the terms in the derivatives of $u$, we obtain part (i).

\vspace{5pt}
\noindent(ii) This follows from the main estimate
(Proposition~\ref{prmest}) and the fact that
$P_\rho$ is obtained by gluing fully elliptic operators.

\noindent(iii) This is deduced by a simple modification of the
arguments used to prove  $f_\rho$ smooth. The details
are omitted.
\end{proof}

\subsubsection{Implicit function theorem} 
\label{ift}Relative to the
decompositions of Proposition~\ref{nlprop}, write $u = (u_1,u_2)$,
$f(0) = (f_1(0),f_2(0))$, $P = (P_{ij})$ and $Q = (Q_1,Q_2)$.  The equation to
be solved becomes the pair
\begin{eqnarray}
f_1(0)  + P_{11}u_1 + P_{12}u_2 + Q_1(u_1,u_2) & = &0,
\label{fdimpart} \\
f_2(0)  + P_{21}u_1 + P_{22}u_2 + Q_2(u_1,u_2) & = &0.
\label{idimpart}
\end{eqnarray}
By the Proposition, $P_{22}$ is invertible, while from construction of
the asymptotic kernels and cokernels in \S\ref{stretch}, the
operator norms of
the other $P_{ij}$ tend to zero as $\rho \to \infty$. Thus for each
fixed $u_1$,
\eqref{idimpart} can be reformulated as a fixed-point problem
$$
u_2 = T_{u_1}(u_2) :=
-P_{22}^{-1}f_2(0)   - P_{22}^{-1}P_{21}u_1 -P_{22}^{-1}Q_2(u_1,u_2).
$$
From Proposition~\ref{nlprop} it is easy to show that if $\rho$ is
sufficiently large and $\|u_1\|$ is sufficiently small, then $T_{u_1}$
is a contraction mapping on a sufficiently small neighbourhood of $0$
in $U''_\rho$. To be more precise there exist $r_1>0$, $r_2>0$ independent of
$\rho$ assumed large, and a function 
$$\varphi_\rho: \{u_1\in U'_\rho: \|u_1\| \leq r_1\} \to
\{u_2\in U''_\rho: \|u_2\| \leq r_2\}
$$
such that every solution $(u_1,u_2)$ of \eqref{idimpart} with
$\|u_j\|\leq r_j$ is of the form $(u_1,\varphi_\rho(u_1))$. 

Thus we have solved an `infinite-dimensional component' of the conformal ASD
equations \eqref{idimpart}; these are reduced to finding zeros of the nonlinear map
between finite dimensional spaces got by substituting $u_2 =
\varphi(u_1)$ into \eqref{fdimpart}: put
$$
\psi_\rho(u_1) = f_1(0)  + P_{11}u_1 + P_{12}\varphi_\rho(u_1) + Q_1(u_1,\varphi(u_1))
$$
and let $\sigma_\rho(u_1)$ be the component of this in $V'_\rho \cap
L^p(E^2)$. To summarize:
\begin{propn} For sufficiently large $\rho$, there exists a nonlinear
map $\sigma_\rho$ from a ball in $U'_\rho$ to $V'_\rho$ whose zeros
correspond to $L^p_2$ conformally ASD metrics near $c_\rho$. In
particular if $V'_\rho  = 0$ such conformally ASD metrics on $X_\rho$
always exist.
\label{p512}\end{propn}

\subsection{Considerations of regularity}

Our final task is to establish that the weak ($L^p_2$) solution found 
in the previous section is actually smooth.  In fact we prove both 
interior regularity and that the resulting metric has optimal boundary 
regularity---in other words that it is polyhomogeneous (relative to an 
index set that we do not specify). In this section $\rho$ is large 
but fixed and we drop it from the notation.

\begin{propn} Let $u\in L^p_2(X)$ satisfy
    \begin{equation} F^w(u) = a,\; w^{-1}L^*(wu) = b
	\end{equation}
where $a$ and $b$ are 	$\ci$,
polyhomogeneous and vanish at the boundary. Then  
if $\sup_X |u|$ is sufficiently small, 
$u\in \ci(X_\rho)$  and is 
polyhomogeneous at $\p X_\rho$.
\label{p513}\end{propn}

\begin{proof}  Combine the equations in the form
\begin{equation}\label{silly}
    P_u u = wQ(u,\nabla^w u) + \theta
\end{equation}
where $\theta$ is a polyhomogeneous section which is a linear
combination of $W^+_0$, $a$ and $b$, 
$$
Q(u,\nabla^w u) = \varepsilon_2(1+wu,\nabla^w u\otimes \nabla^w u) 
+ \varepsilon_3(1+wu, u\otimes u) 
$$
and 
$$
P_u v = P v - wu\varepsilon_1(1+wu,\nabla^w\nabla^w v).
$$

\noindent (i) {\em Interior regularity}.  That $u$ is $\ci$ in $X^o$ 
now follows from standard regularity results, which are applicable 
because $u$ is already $C^{1,\alpha}$, where $\alpha = 1 - 4/p$ 
\cite{ADN}.

\vspace{10pt}
\noindent(ii) {\em Boundary regularity}.  The method we use is closely
analogous to that used by Mazzeo in \cite{RMsy}; we are indebted to him
for useful discussions on this point.

Near the boundary, 
$w = x^\alpha$, 
where $\alpha >0$ and $\theta = x^\alpha\theta'$, say, where
$\theta'$ has a phg asymptotic expansion. Then \eqref{silly} takes the
form
\begin{equation} \label{silly2}
P_u u = P u  - x^\alpha\ve_1(1 + x^\alpha u,\nabla\nabla u)
= x^\alpha Q(u,\nabla u) + x^\alpha \theta'
\end{equation}
and we already know that $\sup(|u| + |\nabla u|)$ is uniformly bounded
as $x \to 0$.   Now if we had the indicial operator $I(P)$ in place of
$P_u$ on the LHS of \eqref{silly2}, we could use the fact that $I(P)$
has an inverse which behaves well on the $b$-Sobolev spaces to
conclude that $u\in x^\alpha L^p_k$ for every $k$ and some fixed
$p$. It would follow that $|(x\p_x)^j\p_y^\beta u|$ is continuous at $\p X$
for all $j$ and all multi-indices $\beta$.  Continuing with the 
argument under the simplifying assumption that $I(P)$ not $P_u$ is on 
the LHS of \eqref{silly2}, we can now
use the fact that $I(P)$ and its inverse also preserve 
spaces of polyhomogeneous functions. To do so, assume by induction that
$u$ has a phg expansion up to some order $N$, say.  Then because of the
factor $x^\alpha$ on the RHS of \eqref{silly2}, $I(P)u$ has an
expansion to order $N + \alpha$, and so $u$ also has such an
expansion. Hence $u$  has a complete phg expansion at the boundary.

The result with $P_u$ replacing $I(P)$ comes from a suitable
approximation argument.  The details of this are straightforward but
lengthy, and are omitted.
\end{proof}

\setcounter{equation}{0}
\section{Main theorems}
\label{mainthm}
\subsection{Gluing ASD $b$-manifolds}

We now summarize our work so far by giving statements of the main
theorems. For the reader's convenience we gather first the relevant
notation.

\subsubsection{Notation} For $j=1,2$, $X_j$ is a 
$4$-manifold with boundary and $g_j$ is a conformally ASD
polyhomogeneous $b$-metric. A piece (union of compact connected
components) $Y$ of $\p X_j$ is given, such that the $g_j$ approach
isometric cylindrical metrics near $Y$.  Weights $w_j$ are chosen as
in \S\ref{wtdef} and a real number $\delta$ is chosen as in \eqref{drange} and
finite-dimensional `cohomology spaces' are defined by the exactness of
the sequence
$$
0 \to H^{1,\pm}_j \to 
w_jx_j^{\pm \delta}L^p_k(X_j,E^1) \to
w_jx_j^{\pm \delta}[L^p_{k-2}(X_j,E^2)\oplus L^p_{k-1}(X_j,E^0)] \to
H^{2,\pm}_j \oplus H^{0,\pm}_j \to 0. 
$$
\label{notn6}
For large real $\rho$, $X_\rho$ is constructed by gluing the $X_j$
across $Y$ and the approximately conformally ASD $b$-metric $g_\rho$
is constructed on $X_\rho$ as in \S\ref{glueb}.

\begin{result}\label{resA}
 Let the notation be as in \S\ref{notn6}. Then for all
sufficiently large $\rho$, there exists a map $\sigma_\rho$ from a
neighbourhood of $0$ in $H^{1,+}_1 \oplus H^{1,-}_2$ into
$H^{2,+}_1\oplus H^{2,-}_2$ such that $\sigma_\rho^{-1}(0)$ parameterizes
the set of conformally ASD metrics $\tilde{g}_\rho$ sufficiently close
to $g_\rho$. Here `sufficiently close' means in particular that 
$$
\sup_{X_\rho}|\tilde{g}_\rho - g_\rho| \to 0 \mbox{ as }\rho \to
\infty
$$
and that $\tilde{g}_\rho$ has a polyhomogeneous expansion near $\p
X_\rho$, such that $\tilde{g}_\rho$ coincides with $g_\rho$ at $\p
X_\rho$.
\end{result}

\subsubsection{Remark} The term `parameterizes' is used in a loose
sense here. The family of conformally ASD metrics given by
$\sigma_\rho^{-1}(0)$ is {\em complete} in the sense that every
diffeomorphism class of conformally ASD metrics sufficiently close to
$g_\rho$ appears in the family. However the gauge action of the
diffeomorphism group has only been fixed up to a finite dimensional
residual gauge freedom and correspondingly the true moduli-space will
in general be got by dividing $\sigma^{-1}_\rho(0)$ by a suitable
compact Lie group.  Further details of this are omitted.

\begin{proof} Combine Propositions~\ref{p512} and \ref{p513}. The
estimate on $\tilde{g}_\rho - g_\rho$ follows from the fact that the
$L^p_2$-norm of $u_2$ in \S\ref{ift}
 is of the same order of
magnitude as $\|W^+_\rho\| = O(e^{-(\eta-\delta)\rho})$ and the
estimate \eqref{sob1}.  The completeness of the family of metrics
constructed comes directly from the implicit function theorem.
\end{proof}

When the boundary components  are spherical space-forms, this Theorem
can be combined with our comparison results to give the following `ASD
desingularization theorem':
\begin{result} Let $(M,g)$ be a compact conformally ASD orbifold, and
let $0\in M$ be a point with a neighbourhood modelled on $\RE^4/a$,
where $a$ is an action of the finite group $\Ga$ with isolated
fixed-point set. If $(N,h)$ is a conformally ASD, WALE space whose
asymptotic region is also modelled by $\RE^4/a$, then there exists a
smooth map $\sigma$ from a neighbourhood of $0$ in 
$H^1_c(M)\oplus H^1_c(N)$ into $ H^2_c(M)\oplus H^2_c(N)$ whose zeroes
give conformally ASD metrics on the connected sum $M\sharp N$ obtained
by joinining the asymptotic region of $N$ to a neighbourhood of $0$ in
$M$. Here the deformation cohomology groups are as in
Theorems~\ref{comp1} and \ref{comp2}.
\label{resB}\end{result}

\begin{proof} Combine Theorem~\ref{resA} with Theorems~\ref{comp1} 
and \ref{comp2}.\end{proof}

Another variant covers the case of a connected sum of compact
conformally ASD orbifolds.

\begin{result} For $j=1,2$, let $(M_j,g_j)$ be a compact conformally ASD orbifold, and
let $0_j\in M_j$ be a point with a neighbourhood modelled on the
origin in $\RE^4/a_j$.  Suppose further that $a_1$ and $a_2$ define
{\em complementary singularities} in the sense that there is an
orientation-reversing linear isometry $\phi$ of $\RE^4$ which intertwines
$a_1$ and $a_2$.  Then there is a
smooth map $\sigma$ from a neighbourhood of $0$ in 
$H^1_c(M_1) \oplus \RE\oplus ({\mathfrak so}_4)^a\oplus H^1_c(M_2)$
into $H^2_c(M_1)\oplus H^2_c(M_2)$ whose zeroes
give conformally ASD metrics on the connected sum $M_1\sharp M_2$ obtained
by joinining at $0_1$ and $0_2$.
\label{resC}\end{result}

\begin{proof}  Take $X_1$ to be the conformal blow up at $\{0_1,0_2\}$
of $M_1 \sqcup M_2$. Take $X_2$ to be the $b$-manifold obtained by
gluing $(-\infty, 0]\times (S^3/a_1)$ to
$(-\infty, 0]\times (S^3/a_2)$ by  identifying $0 \times S^3/a_1$ with
$0 \times S^3/a_2$  by $\phi$. Then $X_2$, viewed as a
$b$-manifold, has boundary equal to $-Y \sqcup Y$, while the boundary
of $X_1$ is $Y \sqcup -Y$.  Applying Theorem~\ref{resA} to
$X_1$ and $X_2$ now gives the conclusion, in view of the calculations
of the deformation cohomology groups for a cylinder given in
\S\ref{pscrv} and Theorem~\ref{comp1}. 
(We have denoted by $a$ the action induced by $a_1$ or
$a_2$ on $S^3\times \RE$.)
\end{proof}

\subsubsection{Remark} The domain of $\sigma_\rho$ in this case has a
natural interpretation, the three summands corresponding to
deformations of the conformally ASD structures of the $M_j$ together
with deformations of the `gluing map' $\phi$. 

\subsubsection{Remark} Clearly Theorems~\ref{resA}, \ref{resB} and
\ref{resC} give existence theorems for conformally ASD metrics if the
indicated obstruction spaces vanish. On the other hand, it is sometimes
possible to calculate the leading term in the map $\sigma$ in terms of
the given data, as in \cite{DF}, and so obtain existence results even
in the presence of obstructions.

\subsection{Gluing hermitian-ASD $b$-manifolds}

The results are very similar to those for ASD conformal structures,
and will follow in the same way (from the methods of \S\ref{s5}) once
it has been explained how to glue hermitian-ASD manifolds. This will
be done in the next few paragraphs. For simplicity we deal only with
the case that the metrics being glued are scalar-flat K\"ahler.

\subsubsection{Notation and assumptions}

For $j=1,2$, let  $(X_j,J_j)$ be complex $b$-manifolds of (real) dimension $4$, and
let $Y\subset \p X_j$ be a piece of the boundary. We make the
assumption that the cylindrical neighbourhoods $U_j$ of $Y$ in $X_j$
(cf.\ \S\ref{data}) are {\em biholomorphic}.   This will be the case,
for example, if $X_1$ is the conformal blow-up of a point in a compact
surface and $X_2$ is obtained by conformal blow-down of $\infty$ in an
asymptotically Euclidean space.

Now introduce $b$-metrics $g_j$ on $X_j$ and assume as in
\S\ref{glueb} that in $U_j$,  $g_j$ approaches a standard cylindrical
metric $g_0$. Assume that the isometry also preserves the complex
structures, so that $g_0$, $g_1$, $g_2$ are all Hermitian with respect
to the given complex structure.

Then we can construct $(X_\rho, g_\rho, J_\rho)$ by gluing just as
before. Note that $J_\rho$ is a genuine integrable complex structure
on $X_\rho$ and that $g_\rho$ is $J_\rho$-hermitian. 

The Riemannian product metric on the cylinder is not necessarily
K\"ahler, but if we assume that $(X_j,c_j,\omega_j)$ is conformally K\"ahler
as in \S\ref{sfk}, so that $\omega_j$ is a 
solution of the twistor equation \eqref{te}, then when we glue, we get
$g_\rho$ as before and $\omega_\rho$ which still defines an integrable
complex structure when rescaled to unit length.  However $\omega_\rho$
is only an approximate solution of the twistor equation defined by the
conformal class of $g_\rho$. We can now introduce weights and repeat 
the work of \S5 for the nonlinear map $F_J$ of 
Proposition~\ref{hasdsum}. This yields the Hermitian analogue of 
Theorem~\ref{resA}:

\begin{result} \label{resD} For $j=1,2$, let $(X_j,g_j,\omega_j)$ be 
conformally Hermitian-ASD $b$ manifolds, where $\omega_j$ satisfies 
the twistor equation defined by $g_j$. Choose weights $w_j$ as in 
\S\ref{wtdef} and define finite-dimensional vector spaces 
$H^{r,\pm}_{j,J}$ by the exactness of the sequence
$$
0 \to H^{1,\pm}_{j,J} 
\to w_jx_j^{\pm\delta}L^p_k(X_j,\Om^{1/4}E^1_J)
\to w_jx_j^{\pm\delta}L^p_{k-2}(X_j,\Om^{3/4}E^2_J)
    \to H^{2,\pm}_{j,J} \to 0.
$$
Suppose further that there is a piece $Y$ of $\p X_j$ such that 
cylindrical neighbourhoods $U_j$ of $Y$ in $X_j$ are biholomorphic.
Then for all sufficiently large $\rho$, there is a map $\sigma_\rho$ 
from a neighbourhood of $0$ in $H^{1,+}_{1,J}\oplus H^{1,-}_{2,J}$
into $H^{2,+}_{1,J}\oplus H^{2,-}_{2,J}$ whose zeroes parameterize 
the conformally hermitian-ASD metrics near $g_\rho$.  Moreover these 
metrics have the same boundary behaviour as the metrics constructed 
in Theorem~\ref{resA}.
\end{result}

Restricting to orbifolds and using the comparison theorem 
(\ref{comp3}), we obtain a result about desingularization (or just 
blow-up) of compact hermitian-ASD orbifolds.

\begin{result} \label{resE}
 Let $(M,g)$ be a compact scalar-flat K\"ahler metric and
let $0\in M$ be a point with a neighbourhood biholomorphic to a 
neighbourhood of $0$ in $\CX^2/a$,
where $a$ is an action of the finite group $\Ga$ with isolated
fixed-point set. If $(N,h)$ is a scalar-flat K\"ahler, WALE space whose
asymptotic region is biholomorphic to a neighbourhood of $\infty$ in
$\CX^2/a$, then there exists a
smooth map $\sigma$ from a neighbourhood of $0$ in 
$H^1_{c,J}(M)\oplus H^1_{c,J}(N)$ into 
$H^2_{c,J}(M)\oplus H^2_{c,J}(N)$ whose zeroes
give scalar-flat K\"ahler metrics
on the connected sum $M\sharp N$ obtained
by joinining the asymptotic region of $N$ to a neighbourhood of $0$ in
$M$. Here the deformation cohomology groups are as in
Theorem~\ref{comp3}.
\end{result} 
    
\begin{proof} It is clear that Theorem~\ref{resD} combined with 
Theorem~\ref{comp3} gives a Hermitian-ASD conformal structure $c$ on 
$M\sharp N$. However it is obvious that $M\sharp N$ is K\"alerian, so 
by \cite{boyer}, there is a K\"ahler representative of the conformal 
class $c$ and this must necessarily be scalar-flat.\end{proof}

\subsection{Weakly asymptotically euclidean scalar-flat K\"ahler
metrics on the blow-up of $\CX^2$}

Finally we give a simple application of Theorem~\ref{resD} which is 
not quite included in Theorem~\ref{resE}.

\begin{result}\label{resF} Let $p_1,\ldots,p_n$ be a collection of 
$n$ distinct points in $\CX^2$. Then the (complex) blow-up $M$ of 
$\CX^2$ at the $p_j$ admits weakly asymptotically euclidean 
scalar-flat K\"ahler metrics.
\end{result}

\begin{proof} Recall first the Burns metric $(B,g_B)$ \cite{exsdm}, 
a weakly asymptotically euclidean scalar-flat K\"ahler metric on the 
blow-up $B$ of $\CX^2$ at the origin. The complement of the 
exceptional divisor in $B$ is biholomorphic to $\CX^2\backslash 0$, 
so we can construct the multiple blow-up $M$, with its standard 
complex structure, by gluing a copy of $B$ at each of the $p_j$. More 
precisely, let $(\check{B},\check{g}_B)$ denote the conformal blow-down of 
$B$. Then $H^2_{c,J}(\check{B}) =0$ by Theorem~\ref{alevanish} and 
Theorem~\ref{comp3}.  So let $X_1$ stand for the $b$-manifold 
obtained from $\CX^2$ by conformal blow-up of each of the $p_j$ and 
conformal blow-down of $\infty$. By Theorem~\ref{comp3},
$$
H^2_{c,J}(X_1) = H^2_{c,J}(\CX^2)\cap O(R^{-2})
$$
and this is zero, because if $Su=0$ in euclidean space then every 
component of $u$ satisfies $\Delta^2 u =0$ and so $u=0$ by the growth 
condition at $\infty$. We take $X_2$ to be $n$ copies of $\check{B}$ 
and apply Theorem~\ref{resD}, taking all weights equal to $1$ since 
all boundary components are spherical (Proposition~\ref{sspec}). The 
conclusion is that $M$ admits a hermitian-ASD conformal structure with 
good asymptotic behaviour at the boundary. It remains to verify that 
there is a weakly asymptotically euclidean K\"ahler representative of 
this conformal class.

The fundamental $2$-form $\omega$ of the metric constructed by 
Theorem~\ref{resD} has the
form $dt\wedge e_1 + e_2\wedge e_3+ O(e^{-t})$, where the $e_j$ form a standard
basis of left-invariant $1$-forms on $S^3$. Following the method of
\cite{boyer}, introduce a $1$-form $\beta$ which measures the failure of
$\omega$ to be K\"ahler: 
$$
d\omega + \beta\wedge\omega = 0.
$$
From the asymptotic formula for $\omega$, $\beta = 2dt + \beta'$ where
$\beta'$ is a $1$-form whose length decays exponentially as $t\to
\infty$.  It is a local calculation that on a Hermitian-ASD manifold
$d\beta$ is ASD hence $\beta$ is harmonic. Because $d\beta = d\beta'$ and
$\beta'$ is in $L^2$, we conclude by the standard integration-by-parts
argument that $d\beta =0$. Hence $\beta'$ is closed and decaying and
so by \cite[Ch.\ 6]{tapsit} represents an element of the 
relative cohomology group
$H^1(M,\p M)=0$. By other results in \cite[Ch.\ 6]{tapsit},
$\beta' = du$ where $du$ is also exponentially decaying.  It follows
that $d(e^{u + 2t}\omega)=0$ and this is the conformal factor that
yields the K\"ahler representative.
\end{proof}

\setcounter{equation}{0}
\section{Existence of ASD conformal structures for manifolds with
boundary}
\label{nct}
Let $(M,g_M)$ be an oriented Riemannian $4$-manifold without 
boundary. Taubes~\cite{taubes-conformal} has shown that there exists 
$N>0$ such that $M_N := M\sharp N\OCP$ admits a conformally ASD metric. 
In this section we shall sketch the adaptations of his argument that 
are needed to prove the analogous result in the $b$-category:
	\begin{result}\label{exist}
Let $X$ be a compact oriented $4$-manifold with
boundary, $g_0$ an exact $b$-metric that is conformally flat near each
component of $\p X$. Then there exists $N>0$ such that $X\sharp N\OCP$
admits a conformally ASD $b$-metric $g$, such that $|g-g_0|\to 0$ at
$\p X$. 
	\end{result}
The strategy of the proof of this result for manifolds without
boundary is summarized in~\cite[Ch.\ 7]{taubes-glueing}. The proof can
be divided into three steps, each of which is substantial:
\vskip 5pt

\noindent\parbox[t]{1.5cm}{Step 1}\hfil
\parbox[t]{13.5cm}{\noindent
Show that for $N>0$ there is a way to construct a Riemannian
metric $g_{N}$ on $M_N$ with the property 
that $\|W^+(g_{N})\| \to 0$ as $N \to \infty$.}
\smallskip

\noindent\parbox[t]{1.5cm}{Step 2}\hfil
\parbox[t]{13.5cm}{Apply the implicit function theorem (IFT) to find a small 
perturbation $g' = g_N1 + h(g_N)$ that is conformally ASD modulo the 
vanishing of a set of constraint functions (essentially the 
map $\sigma$ of 
Proposition~\ref{p512}). Interpret the constraint functions as a finite 
number of nonlinear conditions upon $g_N$.}
\smallskip

\noindent\parbox[t]{1.5cm}{Step 3}\hfil
\parbox[t]{13.5cm}{Show that by replacing $M_{N}$ by $M_{N+n}$, a metric 
$g_{N+n}$ can be constructed with small $\|W^+(g_{N+n})\|$ and vanishing 
constraint functions. Apply the IFT in Step 2 to obtain a 
perturbation of $g_{N+n}$ that is conformally ASD.}
\smallskip

Another way of interpreting Step 3 is to say that one can reduce to 
an unobstructed deformation problem by forming the connected sum with 
sufficiently many copies of $\OCP$, this number depending only upon 
the original data $(M, g_M)$.  An important technical point is that
the norm used in Step~1 is not a standard Sobolev norm; it is
scale-invariant (like $L^2$) but a little stronger. This means
that the IFT used in Step 2 is to be applied in non-standard Banach
spaces; this in turn requires the development of other non-standard estimates for
linear elliptic operators.

\subsection{Sketch Proof I}\label{proof1}
Now let us turn to the $b$-manifold $(X,g_0)$ of Theorem~\ref{exist}. 
We start by applying Taubes's theorem to the double $(M,g_M)$ of $X$. 
That is, regarding $X$ as a non-compact manifold with a cylindrical 
end, we cut off the cylinder and glue the resulting manifold to 
another copy of itself (with opposite orientation). The closed
manifold $M$ then contains a neck $Y\times [-\rho,\rho]_t$ on which the metric 
is a conformally flat Riemannian product metric.  
Here the double is used for definiteness only. Any closed, oriented $4$-manifold 
containing the subset $\{t\leq \rho\}$ of $X$ would do just as well. 

An examination of Step 1 above reveals that $\|W^+\|$ can be decreased 
by gluing copies of $\OCP$ onto $\mbox{Supp}(W^+(g_M))$. In 
particular for our double $M$, we can construct $(M_N,g_N)$ such that 
$\|W^+(g_N)\|\to 0$ as $N\to \infty$, but leaving 
the cylindrical neck in $M$ untouched. 

Now we take Steps 2 and 3 to obtain a conformally ASD metric $g'$, 
say, on $M_{N+n}$. We claim that near the middle ($t=0$) of the 
neck, $g'$ will be a very small perturbation of the product metric. 
Now return to a $b$-manifold $X'$ by cutting the middle of the neck 
and gluing on a semi-infinite cylinder $Y \times [0,\infty)$ to 
produce a manifold with an end $Y\times (-\rho,\infty)$. Glue $g'$ to 
the product metric by means of a cut-off function to obtain a metric 
$g''$ on $X'$ with the property that $W^+(g'')\not=0$ only in a small 
neighbourhood of $t=0$. We can moreover assume that a suitable 
weighted Sobolev norm of  $W^+(g'')$ is  as small as we please.

Thus we are now in the same framework as for the deformation
theory in the rest of this paper. So we invoke once more 
the IFT to find a conformally ASD $b$-metric $g'''$ 
as a small perturbation of $g''$. We
claim that arguments analogous to those of Step~3 allow us to overcome
the obstructions that could arise here.

The inelegant double use of the IFT in this argument is intended to
avoid the need to adapt to $b$-manifolds the non-standard norms
mentioned above. This  completes our sketch of a proof 
of Theorem~\ref{exist}.

\subsection{Sketch Proof II}\label{proof2}
It is also possible to argue slightly differently: apply Step~1 as 
outlined above, and then pass back to a $b$-manifold $X'$ by cutting at 
$t=0$ and gluing in $Y\times [0,\infty)$. Now adapt Steps~2 and 3 to 
apply to $b$-manifolds with small $\|W^+\|$. This direct 
approach is attractive conceptually and several of the main steps go 
through without major changes. However, it is technically subtle  as
we have already indicated because of the non-standard norms used
throughout Taubes's argument.  In particular, Taubes's  analysis 
makes
heavy use of the spectral theory  
of (bundle-valued) Laplacians $\nabla^*\nabla$: for some estimates
it is necessary to expand sections as a linear combination of 
eigensections of $\nabla^*\nabla$, and spectral projection is used 
to define  finite-dimensional subspaces corresponding to
`small eigenvalues'.  Because $\nabla^*\nabla$ 
has continuous spectrum on a $b$-manifold, it is not simple to extend 
such arguments to $b$-manifolds. However we claim that one can
satisfactorily glue Taubes's 
estimates over the compact piece $t\leq 0$ of $X'$ onto standard
weighted Sobolev space estimates for a product metric on the
cylindrical end $Y\times (0,\infty)$. This gives another approach to
the proof of Theorem~\ref{exist}.

\setcounter{equation}{0}
\section{Vanishing theorems}
\label{vanishing}
This section is summarizes some vanishing theorems for the
obstruction-spaces that arise in the main gluing theorems. They apply
to conformally ASD $b$-manifolds that are either obtained by conformal
blow-up or blow-down of a conformally ASD manifold whose metric has
additional geometric properties, for example, Einstein or scalar-flat K\"ahler.

The analysis of the Einstein cases
rests on the use of Weitzenbock formulae for $D^*$, while the analysis
of the scalar-flat K\"ahler story was done in~\cite{LebS} in the
compact case and is a modification of this in the ALE case.

\subsection{$D^*$ and Dirac operators}

In \S\ref{ASDdef} we have described $D$ in terms of a coupled version of 
the operator $d^+d^*$. There is a useful alternative 
description using spinors, which we now explain.
Spinors cannot be introduced globally on a $4$-manifold unless the
second Stiefel--Whitney class vanishes, but we shall only ever need
tensor products of an even number of spin-bundles, and these always
exist globally. Indeed the link between our two accounts of the
operator $D_g$ is obtained precisely by carrying through the necessary
identifications of tensor products of spin-bundles with certain
bundles of tensors (associated to the tangent bundle) \cite{OT}.

On an oriented Riemannian $4$-manifold, the spin-bundles $V^{\pm}$ are
naturally $SU(2)$-bundles and that the Dirac operators interchange $+$ and $-$:
$$
\dir^+:C^\infty(V^-) \to C^\infty(V^+), \;\;
\dir^-:C^\infty(V^+) \to C^\infty(V^-).
$$
Again there are coupled versions of these,
$\dir^\pm_E:C^\infty(V^{\mp}\otimes E) \to C^\infty(V^{\pm}\otimes
E)$, for any vector bundle $E$ equipped with a unitary connection. 
In particular, there is a second-order operator $C^\infty(S^2 V^-) \to
S^2(V^+)$ given by composing the Dirac operators $C^\infty(V^-\otimes
V^-)\to C^\infty(V^+\otimes V^-)$ and $C^\infty(V^+\otimes V^-) \to
C^\infty(V^+\otimes V^+)$. Up to a constant factor, this can be
identified with $d^+\delta$, using the natural isomorphisms $S^2V^{\pm}
= \Lambda^{\pm}$. 
Now ignoring conformal weights,
$$
E^0 = V^+ \otimes V^-,\;\;\;
E^1 = S^2V^+ \otimes S^2V^-,\;\;\;
E^2 = S^4V^+.
$$
(Here the canonical symplectic forms on $V^{\pm}$ have been used to
eliminate all appearances of dual spin spaces.)

Then $D^*_g$ can be written in terms of coupled Dirac operators $D_1$
and $D_2$ where
$$
D_1:C^\infty(V^-\otimes S^3V^+) \to C^\infty(V^-\otimes V^-\otimes
S^2V^+)
$$
and
$$
D_2:C^\infty( S^4V^+) \to C^\infty(V^-\otimes S^3V^+);
$$
namely $D^*_g = S(D_1D_2) + \Phi$.
The operator $S$ is the algebraic operation of symmetrization:
$V^-\otimes V^-\otimes S^2V^+ \to S^2V^-\otimes S^2V^+$.  There is a
similar formula for $D_g$, which we do not write down. 

The above formula for $D^*_g$ has an important simplification if $g$
is ASD. This is that the symmetrization is {\em unnecessary}. The
reason is that the skew part of $D_1D_2$ is an algebraic
operator $S^4V^+ \to \Lambda^2 V^-\otimes S^2 V^+ = S^2V^+$, given by 
partial contraction with some component of the curvature tensor. In
general,
the only component that could provide such a map  is $W^+ \in
C^\infty(S^4V^+)$, which we have assumed is zero. The conclusion is as
follows: if $g$ is ASD, then $D_g^*$ can be identified with the
operator 
$$
D_1D_2 + \Phi : C^\infty(S^4V^+) \to C^\infty(S^2V^-\otimes S^2V^+).
$$
\subsection{Vanishing theorems when $X$ is the conformal
blow-up of a compact ASD-Einstein orbifold}

To simplify notation, let us drop the bars which have been used to
distinguish a  compact manifold from its conformal blow-up; this
should cause no confusion since the latter will not be used in this section.
Our first result is the following:
\begin{propn}
Let $(X,g)$ be a compact $4$-orbifold such that $g$ is ASD and
Einstein with positive scalar curvature $s$. Then $H^2_c(X) =0$. If
instead $(X,g)$ is ASD and Ricci-flat, then $H^2_c(X)$ consists of
parallel sections of $E^2$.
\end{propn}
This is a folklore theorem, but we reproduce the short proof.
	\begin{proof}
 If $X$ is Einstein, we can identify $D_g^*$ with
the composite $D_1D_2$ of Dirac operators as above.
To prove the proposition, it is enough to note the Weitzenbock
formulae 
$$
D_1^*D_1 = \nabla^*\nabla + \frac{5}{12}s ,\;\;
D_2^*D_2 =\nabla^*\nabla + \frac{1}{2}s
$$
which hold whenever $\Phi=0$ and $W^+ =0$.  (The verification of these
is left to the reader.) Suppose $D_g^*\psi=0$; set $\chi =
D_2\psi$, so that $D_1\chi=0$. If $s>0$ then  $D_1^*D_1$ is
invertible, so $\chi=0$, i.e.\ $D_2\psi=0$. Similarly $D_2^*D_2$ is
invertible, so $\psi=0$. This completes the proof of the first part.

Suppose now that $s=0$. With $\chi$ and $\psi$ as before, we
deduce first that $\chi$ is a {\em parallel} section. In particular,
$D_2^*\chi =0$. But $\chi$ also lies in the image of $D_2$ (by
definition), so $\chi$ must be zero, by the Fredholm alternative.
Thus $D_2\psi=0$ and, applying the Weitzenbock formula, 
$\psi$ is
parallel. Thus we have identified the kernel of $D_g^*$ with the space
of parallel sections of $S^4V^+$ in the Ricci-flat case.
	\end{proof}

Of course compact ASD-Einstein manifolds with $s \geq 0$ are rather
rare. A well-known result of Hitchin states that the only examples
with $s>0$ are the complex projective plane (with the opposite
orientation) and the $4$-sphere. And when $s=0$, one has only
(quotients of) the K3-surface or the flat $4$-torus.  Since the former
is simply connected and the latter has trivial holonomy, we have $\dim
H^2_c = 5$ in each of these cases.

When the class of spaces is widened to compact, ASD--Einstein {\em
orbifolds} with $s>0$, many more examples appear. These include the
weighted projective spaces of Galicki--Lawson \cite{GL}.

The next class of examples consists of  the compact scalar-flat K\"ahler
surfaces. Here there is a fundamental dichotomy according as the Ricci
tensor does or does not vanish.  Since the Ricci-flat case was already
analyzed, we may as 
well suppose that the surface is not Ricci-flat. Then we have the
vanishing theorem of \cite{KLP} to the effect that $H^2_c(X)=0$
whenever the scalar-flat K\"ahler surface is {\em non-minimal}; that
is to say, whenever it contains at least one rational curve of
self-intersection $-1$. 
We shall not repeat the argument
(though the reader will see many of the details in our discussion
below of the case of non-compact but ALE scalar-flat K\"ahler
surfaces). 

\begin{remark}
 Although we do not give a formal statement, it also often
possible to compute $H^2_c$ when $c$ is conformally flat. Indeed, in
this case the ASD deformation theory is essentially the same as the
conformally flat deformation theory; in particular if the latter is
unobstructed, then so is the former. On the other hand, the
conformally flat deformation theory is given by a flat complex (de
Rham complex with twisted coefficients). The cohomology groups can
sometimes be computed by topological methods. For example, for 
generic conformally flat structures, 
$H^2_c(S^1\times S^3\sharp\cdots\sharp S^1\times S^3)=0$ 
by a Meyer--Vietoris argument.
\end{remark}

\subsection{Vanishing theorems for WALE spaces}

In this subsection we consider the case that the cylindrical-end model
has a conformal blow-up that is WALE and either Ricci-flat or
scalar-flat K\"ahler. Once again, we shall not have any use for the
$b$-manifold here and therefore drop the use of hats.

Recall from previous discussion that we are
interested in the part of the kernel of $D_g^*$ that is $O(R^{-2})$
near $\infty$.

\begin{propn} Suppose $(X,g)$ an ALE space that is ASD and
Ricci-flat. Then $H^2_g(X)=0$.
\end{propn}
	\begin{proof}
 Recall the formula $D^* = D_1D_2$ from the
previous proposition. Let $\psi \in H^2_g$ and let $\chi =
D_2\psi$. Then $|\chi| = O(R^{-3})$ as $R\to \infty$, and so
$$
0=\int_{R \leq R_0}(\chi,D_1^*D_1\chi) =
\int_{R \leq R_0}(\chi,\nabla^*\nabla\chi)
= \int_{R \leq R_0} |\nabla\chi|^2 + O(R_0^{-4})
$$
and letting $R_0 \to \infty$ we conclude as
before that $\chi$ is parallel. Since $|\chi| \to 0$ at $\infty$,
moreover,  $\chi=0$. Hence $D_2\psi=0$. We make the same argument as
above using the Weitzenbock formula. This time the boundary term is
$O(R_0^{-2})$ but we still conclude that $\psi$ is parallel and hence
$0$, since $0$ at $\infty$.
	\end{proof}

Now let us take up the K\"ahler story. By definition a non-compact 
K\"ahler surface $(X,g,J)$ is said to be ALE if $g$ is ALE and
K\"ahler with respect to $J$ {\em and if} the chart at infinity $\phi$
can be chosen to be a biholomorphic map 
$X-K \to (\CX^2 - B)/\Ga$. We shall now show that in this case
$H^2_g(X)=0$. Our method is to
analyze $D^*_g$ in
terms of holomorphic data on $X$, decaying  at $\infty$. This is
helpful because of the following lemma:

\begin{lemma} Suppose that $X$ is an ALE K\"ahler surface and suppose
that $T$ is a holomorphic tensor field on $X$, $|T| \to 0$ at
$\infty$. Then $T=0$.
\label{holv}
\end{lemma} 
\begin{proof}
 Transfer $T$ to $(\CX^2-B)/\Ga$ and pull back to
$\CX^2-B$. Then $T$ becomes a holomorphic section of a trivial vector bundle
and so each component of $T$ is a holomorphic function  that decays at
$\infty$. But by the removable singularity theorem of Hartogs, each
component extends uniquely to a holomorphic function on $\CX^2$ and by
the maximum principle must therefore be identically zero. This
argument shows that $T$ is identically zero on $X - K$. But then by
uniqueness of analytic continuation, $T$ is identically zero on
$X$.
	\end{proof}

Now the main theorem of this subsection can be given:

\begin{thm}  Suppose that $(X,g,J)$ is an ALE scalar-flat K\"ahler
surface. Then $H^2_g(X)=  H^2_{g,J}(X)= 0$.
\label{alevanish}\end{thm}
	\begin{proof}
Recall first that
$$
\Lambda^+\otimes\CX = K \oplus 1 \oplus K^{-1}, \;\;
\Lambda^-\otimes\CX  = \Lambda^{(1,1)}_0.
$$
Here $K$ is the canonical-bundle of $X$, and the trivial bundle $1$ is
embedded in $\Lambda^+$ by multiplication by $\omega$. Let us for the
moment write $K\oplus K^{-1}=E$ so that $\Lambda^+ = 1 \oplus E$. This
decomposition is preserved by the connection since $J$ is parallel,
so that $D^*$ decomposes as a pair of
operators
\begin{equation} \label{op1}
S^*= d^-\delta + \Phi : \ci(\Lambda^+)\to \ci(\Lambda^-)
\end{equation}
(taking the component of $\omega$) and 
\begin{equation} \label{op2}
d^-_E\delta_E + \Phi : \Omega^+(E) \to \Omega^-(E).
\end{equation}
Now we shall decompose the other factor of $\Lambda^+$. If $\psi$ is a
section of $E^2$, write its components as follows
$$
\begin{pmatrix}
\psi_0 & \psi_{1}\\
         \psi_{1} & \psi_{11}
\end{pmatrix} \in C^\infty 
\begin{pmatrix}
1 & 1\otimes E \\
 E\otimes 1 & E\otimes E
\end{pmatrix}.
$$
Now the condition $D_g^*\psi=0$ may be written as the two equations
\begin{equation} \label{st6}
d_-d^c \psi_0 + \rho \psi_0 + d_-\delta \psi_1 = 0
\end{equation}
and 
\begin{equation} \label{st7}
d_-^E(d^E)^c \psi_1 + \rho \psi_1 + d_-^E\delta^E \psi_{11} = 0
\end{equation}
The vanishing theorem will be proved according to the following scheme: 
\begin{quote}
(i) $H^0(X,\Theta)\cap O(R^{-2}) =0$ implies $\psi_0 =0$;

(ii) $H^0(X,\CO(K))\cap O(R^{-2}) =0$ implies $\psi_1 = 0$;

(iii) $H^0(X,\Omega^1\otimes K)\cap O(R^{-2}) =0$ and
$H^0(X,\CO( K^2))\cap O(R^{-2}) =0$  implies $\psi_{11} =0$.
\end{quote}
In other words, relative to $S^2_0\Lambda^+ = K^{-2}\oplus K^{-1} \oplus
1 \oplus K \oplus K^2$, we eliminate first the component in the
trivial bundle, next the components in $K^{\pm 1}$, finally those in
$K^{\pm 2}$. Remark that complex conjugation carries $K^r$ into
$K^{-r}$ so it is enough to deal with the components in $1$, $K$ and
$K^2$.

\vspace{10pt}
\noindent {\em Proof of (i)} As in~\cite{LebS}, the real function
$\psi_0$ satisfies 
Lichnerowicz's differential equation; on a compact manifold it follows
that $\nabla^{1,0}\psi_0$ is a holomorphic
vector field. The argument requires integration by parts but in our
situation we have sufficient decay at $\infty$ so that the conclusion
holds. In fact, $\nabla^{1,0}\psi_0$ is  holomorphic and 
decays at $\infty$, so by Lemma~\ref{holv} $\psi_0$ is
constant. Finally  $\psi_0$ is
$O(R^{-2})$ so $\psi_0=0$; the proof of (i) is complete. 

\vspace{10pt}
\noindent {\em Proof of (ii)}.  Referring to (\ref{st6}), $\psi_1$ satisfies the
equation $d^-\delta \psi_1 =0$.  As in the compact case,
this implies $d\delta\psi_1=0$ (cf.\ the proof of the existence of the
conformal factor). In particular $\psi_1$ is harmonic.
Because the Hodge and $\db$-Laplacians agree (up to a
factor of $2$) we infer that $\db^*\db\gamma =0$
where $\psi_1 = \gamma + \overline{\gamma}$ is the decomposition of
according to components in $K$ and $K^{-1}$. Integration-by-parts is
applicable now to show that $\gamma$ is holomorphic. The proof is
completed by applying Lemma~\ref{holv} to $\gamma$.

Finally we consider implication (iii). Since $\psi_0=0$ and
$\psi_1=0$, $\psi_{11}= \alpha + \overline{\alpha}$, say, where $\alpha$
is a section of $K^2$, and (\ref{st7}) gives $d^-_K\delta_K\alpha
=0$. These operators are the usual Hodge-de Rham operators, coupled
to the holomorphic line bundle $K$ and  $\alpha$ lies in the space
$\Omega^{2,0}(K)$. For reasons of degree, and using the K\"ahler
identity $\partial^* = i[\Lambda,\db]$,
$$
d_Kd_K^*\alpha = id_K\Lambda\db_K\alpha = \omega\lambda + \mu
$$
where $\lambda$ is a section of $K$ and $\mu$ is a section of
$K^2$. More precisely,
\begin{equation} \label{st10}
i\db_K\Lambda\db_K\alpha = \omega\lambda.
\end{equation}
Since $\db_K^2=0$ (the curvature is of type $(1,1)$), we obtain
$\omega\wedge\db_K\lambda =0$ and hence $\lambda$ is a holomorphic
section of $K$. Since $\lambda$ decays at $\infty$, we have $\lambda
=0$. Two more steps, each involving an application of Lemma~\ref{holv}
complete the proof. For with $\lambda=0$,
equation (\ref{st10}) says that 
$\Lambda\db_K\alpha$ is a decaying holomorphic section of
$\Lambda^{1,0}(K)$, hence zero. Since $\Lambda:\Omega^{2,1} \to
\Omega^{1,0}$ is an isomorphism, it follows that $\alpha$ is a
decaying, holomorphic section of $K^2$. This completes the proof
of (iii) and  
hence the vanishing theorem.
	\end{proof}

\end{document}